\documentclass[a4paper,12pt,twoside, fleqn]{amsart}
\usepackage{amsmath,amsthm,amsfonts}
\usepackage[T1]{fontenc}
\newtheorem{hyp}{Hypothesis}
\newtheorem{thm}{Theorem}[section]

\newtheorem{cor}[thm]{Corollary}
\newtheorem{lem}[thm]{Lemma}

\newtheorem{prop}[thm]{Proposition}
\newtheorem{proposition}[thm]{Proposition}
\newtheorem{definition}[thm]{Definition}

\newtheorem{defi}[thm]{Definition}

\newtheorem{nota}[thm]{Notation}

\theoremstyle{definition}

\newtheorem{rem}[thm]{Remark}
\newtheorem{rems}[thm]{Remarks}

\newcommand{\R}{{\mathbb{R}}}

\newcommand{\E}{{\mathbb{E}}}

\newcommand{\C}{{\mathbb{C}}}

\newcommand{\F}{{\mathbb{F}}}
\newcommand{\Z}{{\mathbb{Z}}}

\newcommand{\T}{{\mathbb{T}}}
\newcommand{\PP}{{\mathbb{P}}}

\newcommand{\Cal}{\mathcal}
\newcommand{\cal}{\mathcal}

\def \Var {{\rm Var}}

\def\cf{{cf. \hskip -2mm}}

\def\j{{\underline j}} \def\k{{\underline k}}

 \def\el{{\underline \ell}}
\def\t{{\underline t}} \def\n {{\underline n}}
\def\0{{\underline 0}}  \def\1{{\underline 1}}

 \def\p {{\underline p}}
\def\u{{\underline u}} \def\r{{\underline r}}
\def\v{{\underline v}} 
\def\s{{\underline s}} 
 \def\q{{\underline q}}

\def\al {{\underline \alpha}}

\def \and{\text{ and }}

\def \stm0{{\setminus \{\0\}}}

\def \eop{\qed}
\def \Proof {\vskip -3mm {{\it Proof}. }}
\def \proof {\vskip -3mm {{\it Proof}. }}
\def \Proof {\vskip -3mm {{\it Proof}. \ }}
\def \Var  {{\rm Var }}

\def \Card {{\rm Card}}

\begin{document}

\date{\today}

\parskip=2mm
\baselineskip 15pt
\parindent=0mm

\title[On the quenched  functional CLT in random sceneries] {On the quenched functional CLT \\ in random sceneries}
\vskip 2mm

\author{Guy Cohen and Jean-Pierre Conze}
\address{Guy Cohen, \hfill \break Dept. of Electrical Engineering, \hfill \break Ben-Gurion University, Israel} \email{guycohen@bgu.ac.il}
\address{Jean-Pierre Conze, \hfill \break IRMAR, CNRS UMR 6625, \hfill \break University of Rennes I, Campus de Beaulieu, 35042 Rennes Cedex, France} \email{conze@univ-rennes1.fr}
\subjclass[2010]{Primary: 60F05, 28D05, 22D40, 60G50; Secondary: 47B15, 37A25, 37A30} \keywords{quenched functional central limit theorem, $\Z^d$-action, random walk in random scenery, self-intersections of a r.w., toral automorphisms,  exponential mixing, flows on homogeneous spaces, $S$-unit, cumulant}

\maketitle

\vskip 3mm 
\begin{abstract}  
We prove a quenched functional central limit theorem (quenched FCLT) for the sums of a random field (r.f.) along a $\Z^d$-random walk in different frameworks: 
probabilistic (when the r.f. is i.i.d. or a moving average of i.i.d. random variables) and algebraic 
(when the r.f. is generated by commuting automorphisms of a torus or by commuting hyperbolic flows on homogeneous spaces).
\end{abstract}

\tableofcontents

\section*{\bf Introduction} 

Let $G$ be a group acting on a probability space $(E, \cal B, \mu)$ by a measure preserving action $(g, x) \in G \times E \to T^g x \in E$.
A random walk $(Z_n)$ defined on a probability space $(\Omega, \PP)$ with values in the group $G$ induces a random walk on $E$.
For $f \in L^2(E, \cal B, \mu)$, we can consider the sums along the random walk: $\sum_{k= 0}^{n-1} f( T^{Z_k} x)$. 

This general framework leads in practice to different situations and methods in the proof of a central limit theorem (CLT) and a functional central limit theorem (FCLT) 
along the paths of the random walk. In particular the proof of the tightness for the FCLT requires specific tools which it seems interesting to present in examples.

A first situation is that of a random walk in random sceneries (cf. \cite{KeSp79}, \cite{Bo89}).
For $d \geq 1$, let $X = (X_\el)_{\el \in \Z^d}$ be a strictly $\Z^d$-stationary real random field (r.f.). 
One can assume that the r.v.s $X_\el$ are defined on a probability space $(E, \cal B, \mu)$ on which commuting measure preserving maps $T_1, ..., T_d$ act in such a way
\footnote{Underlined letters represent elements of $\Z^d$ or $\T^d$. We write $\el$ for $(\ell_1, ..., \ell_d)$ and $T^\el$ for $T_1^{\ell_1}... T_d^{\ell_d}$.
The euclidean norm of $\el \in \Z^d$ is denoted by $|\el|$ or $\|\el\|$.}
that $X_\el = T^\el X_\0$.  

Conversely, given commuting measure preserving invertible maps $T_1, ..., T_d$ and a measurable $f$ on a probability space $(E, \cal B, \mu)$,
$(T^\el f)_{\el \in \Z^d}$ is a strictly $\Z^d$-stationary random field.
If $(Z_n)$ is a random walk in the group $\Z^d$, then the sums along $Z_n$ read $\sum_{k= 0}^{n-1} X_{Z_k} = \sum_{k= 0}^{n-1} T^{Z_k}X_\0$,
or $\sum_{k= 0}^{n-1} T^{Z_k} f$.
When $(X_\el)$ is a $d$-dimensional random field of i.i.d. random variables, we obtain the classical random walk in random sceneries.

Another types of examples in the algebraic case can be obtained as follows.
Suppose that $G = \rm{SL}(\rho, \Z)$ and that $(E, \cal B, \mu)$ is the torus $\T^\rho$, $\rho \geq 2$, endowed  with the Borel $\sigma$-algebra and the Lebesgue measure.
The map $x \to Ax$, where $A$ is a matrix in $\rm{SL}(\rho, \Z)$, defines an automorphism of $\T^\rho$ which preserves $\mu$.
When a spectral gap property is available for the transition operator associated to the random walk on $G$, 
the previous sums for $f$ in a convenient class of observables satisfy a CLT (cf.\cite{AyLiSt09}, \cite{CoLe9}) for $\PP$-a.e. $\omega \in \Omega$.

The commutative case, which we will consider here, is different from the spectral point of view.
For the action of commuting matrices in $\rm{SL}(\rho, \Z)$ acting on  $\T^\rho$, we prove  for $\PP$-a.e. $\omega$ a  functional CLT, 
extending previous results in \cite{CohCo17}.
A second algebraic example comes from commuting flows on homogeneous spaces. 
Based on the exponential mixing of all orders proved in \cite{BjEiGo17}, a CLT has been shown in \cite{BjGo20} for ergodic sums on {F\o{}lner sets 
when the observables are smooth. Likewise we prove here a CLT and its functional version for the sums along  a random walk.

The result, a functional CLT for the different models described above, is presented for a general aperiodic random walk in dimension $d > 1$ 
with a moment of order 2,  but the detailed proofs are given in the case of a centered  2-dimensional r.w.
The proofs can be adapted easily to the case of transient random walks.
We say also some words in the i.i.d. case, when the usual random walk is replaced by a plane Lorentz process 
generated by a periodic billiard with dispersive obstacles (cf. \cite{Pen09}, \cite{Pen14}). 

Beyond the CLT, tightness is a main step in the proof of a FCLT. 
To show it we use the method based on the maximal inequality for associated r.v.s due to Newman and Wright \cite{NeWr81} 
or, in the algebraic case, the method based on norm estimates for the maximum of partial sums (cf. Billingsley \cite{Bi99}, M\'oricz \cite{Mo76}).
A difficulty which occurs is that the estimates available for the random walk involve constants depending on the trajectory.

The content of the paper is the following. 
Section \ref{varGene1} contains results on the variance of sums along a random walk. The independent case is presented in Section \ref{indepSec}.
Some facts on cumulants are recalled in Section \ref{Sectcumulant}, then applied to moving averages in Section \ref{movAverag} and to the algebraic models 
in Sections \ref{algMod1} and \ref{expMixSect}.  For the tightness in the latter cases, we use the method of maximum of partial sums 
in an adapted version presented in Section \ref{SectMoricz}.

The paper is an extension of a previous version. We have added the FCLT along a random walk for flows on homogeneous spaces,
using the recent results in \cite{BjEiGo17} and \cite {BjGo20} on the multiple mixing and on the CLT for group actions which are exponentially mixing of all orders.
We have also added some remarks about the non nullity of the variance, in particular, the observation that there is no degeneracy for the sums along a transient random walk.

\vskip 3mm 
{\bf Acknowledgements.} This research started during visits of the first author to the IRMAR at the University of Rennes 1 
and of the second author to the Center for Advanced Studies in Mathematics at Ben Gurion University. 
The authors are grateful to their hosts for their support.
They thank Y. Guivarc'h and B. Bekka for helpful discussions.

\section{\bf Summation along a r.w. and variance} \label{varGene1}

\subsection{\bf Random walks and sums along random walks} \label{randWalk}

\

First we recall some definitions and results about the random walks on $\Z^d$ (see \cite{Sp64}, details on the results recalled here can also be found in \cite{CohCo17}). 

Let $(\zeta_i)_{i \geq 0}$ be a sequence of i.i.d. random vectors on a probability space $(\Omega, \, \PP)$ with values in $\Z^d$ and common probability distribution $\nu$.
The associated {\it random walk} (r.w.) $Z = (Z_n)$ in $\Z^d$ starting from $\0$ is defined by $Z_0 := \0$, $Z_n := \zeta_0 +... + \zeta_{n-1}$, $n \geq 1$.

The r.v.s $\zeta_i$ can be viewed as the coordinate maps on $(\Omega, \, \PP)$ obtained as  $(\Z^d)^\Z$ equipped with the product measure $\nu^{\otimes \Z}$
and with the shift $\theta$ acting on the coordinates. 
We have $\zeta_i = \zeta_0 \circ \theta^i$ and the cocycle relation $Z_{n +n'} = Z_n + Z_{n'} \circ \theta^n, \forall n, n' \geq 0$.

Let denote by $\cal S := \{\el \in \Z^d: \PP(\zeta_0 = \el) > 0\}$ the support of $\nu$ and by $L$ the sub-lattice of $\Z^d$ generated by $\Cal S$.
Without loss of generality, we can assume that $Z$ is reduced which means that $L$ is cofinite in $\Z^d$.
Therefore the vector space generated by $L$ is $\R^d$ and $d$ is the `genuine' dimension of the random walk $Z$.

For simplicity, we will assume that $L = \Z^d$ (the random walk $Z$ is said to be {\it aperiodic}). 
Observe that one can replace a reduced r.w. $Z$ by an aperiodic one, again without loss of generality.

Let $D$ be the sub-lattice of $\Z^d$ generated by $\{\el - \el', \el, \el' \in \cal S\}$.
We denote by $\Gamma_1$ the annulator in $\T^d$ of $D$, that is the closed subgroup of $ \T^d$ defined by 
$\{\t \in \T^d: e^{2\pi i \langle \r, \t \rangle} = 1, \forall \r\ \in D \}$ and by $d\gamma_1$ the Haar probability measure of the group $\Gamma_1$.
The r.w. is said to be {\it strictly (or strongly) aperiodic}, if $D = \Z^d$.

\vskip 3mm
\goodbreak
{\bf Sums along random walks} \label{sumsAlong}

Given a strictly $\Z^d$-stationary random field $X = (X_\el)_ {\el \in \Z^d}$, where the real random variables $X_\el$ are defined on a probability space $(E, \cal B, \mu)$,
the process of `ergodic sums' along the random walk $(Z_n)$ is
\begin{eqnarray}
&&S_n^{\omega, X}(x) = S_n^{\omega}(x) :=\sum_{k=0}^{n-1} X_{Z_k(\omega)}(x), \, n \geq 1, \, \omega \in \Omega. \label{defsumRW}
\end{eqnarray}
If the random field is represented as $X_\el = T^\el f$, where $T_1, ..., T_d$ are commuting measure preserving maps and $f \in L^2(E, \cal B, \mu)$, the sums read:
\begin{eqnarray}
S_n^\omega f = \sum_{k=0}^{n-1} T^{Z_k(\omega)} \, f = \sum_{\el \in \Z^d} w_n(\omega, \el) \, T^{\el} \, f, \label{defwn0}
\end{eqnarray}
where $w_n(\omega, \el)$ (denoted also by  $w_n^\omega(\el)$) is the local time of the random walk at time $n \geq 1$:
\begin{eqnarray}
w_n(\omega, \el) &=& \#\{k < n: Z_k(\omega) = \el\} =  \sum_{k=0}^{n-1} 1_{Z_k(\omega) =\el}. \label{loctim}
\end{eqnarray}

Summing along the random walk amounts to take the ergodic sums for the skew product
$(\omega, x) \to T_{\zeta_0}(\omega, x) = (\theta \omega, T^{\zeta_0(\omega)} x)$ on $\Omega \times E$. Putting 
$F_f(\omega, x) = F(\omega, x) = f(x)$ for an observable $f$ on $E$, we get that the ergodic sums of $F$ for $T_{\zeta_0}$ read:
\begin{eqnarray}
S_n F (\omega, x) = \sum_{k=0}^{n-1} F(T_{\zeta_0}^k \, (\omega, x)) =  \sum_{k=0}^{n-1} f(T^{Z_k(\omega)} x) = (S_n^\omega f)(x). \label{ergSum1}
\end{eqnarray}
A limit theorem in distribution for the sums $S_n^\omega f$ (with respect to the measure $\mu$ on $E$) obtained for $\PP$-a.e. $\omega$ 
is sometimes called {\it quenched}. We will use this terminology
\footnote{We follow here the terminology of \cite{AyLiSt09} used in several papers.
The term `quenched' is also used in the random scenery when a limit theorem is shown for the distribution with respect to $\omega$, conditionally to the scenery $X$.}.
If the random variables $S_n F (\omega, x)$ are viewed as defined on $\Omega \times E$ endowed with the probability $\PP \times \mu$, 
a limit theorem under $\PP \times \mu$ for theses sums is called {\it annealed}.

\vskip 2mm
\subsection{\bf Variance for quenched processes} \label{varRWSect1}
\ 

Let $f$ be a function in $L^2(E, \cal B, \mu)$ with real values. Everywhere we assume (or prove) the absolute summability of the series of decorrelations
\begin{eqnarray}
\sum_{\el \in \Z^d} | \int_X T^\el f \, \overline f \, d\mu| < \infty, \label{sumcorr1}
\end{eqnarray}
which implies existence and continuity of the spectral density, the even function given by
\begin{eqnarray}
\varphi_f(\t) = \sum_{\el \in \Z^d} \langle T^\el f \, f \rangle \, e^{2\pi i \langle \el, \t \rangle}. \label{specdens1}
\end{eqnarray}
The computation of the variance $\int_E |\sum_{k=0}^{n-1} \, T^{Z_k(\omega)} f|^2 \, d\mu$, is related to the 
number of self-intersections of the random walk at time $n \geq 1$:
\begin{eqnarray}
&&V_n(\omega) := \#\{0 \leq u, v < n: \ Z_u(\omega) = Z_v(\omega)\} = \sum_{\el \in \Z^2} w_n(\omega, \el)^2 
=\int_{\T^d} |\sum_{\el \in \Z^d} w_n(\el) \, e^{2\pi i \langle \el, \t \rangle}|^2 \, d\t. \label{selfinter}
\end{eqnarray}
Let us consider the kernels (which are even functions)
\begin{eqnarray}
&&K(w_n^\omega)(\t) = |\sum_{k=0}^{n-1} e^{2\pi i \langle Z_k(\omega), \t \rangle} |^2 = |\sum_{\el \in \Z^d} w_n(\omega, \el)
e^{2\pi i \langle \el, t \rangle}|^2, \, \tilde K(w_n^\omega)(\t) = V_n(\omega)^{-1} K(w_n^\omega)(\t). \label{kernOmNorm}
\end{eqnarray}
We say that the summation along the r.w. $Z$ is $\xi$-{\it regular}, where $\xi$ is a probability measure on $\T^d$, if (for $\PP$-a.e. $\omega$) 
the normalised kernel $(\tilde K(w_n^\omega))_{n \geq 1}$ converges weakly to $\xi$, 
i.e., $\lim_{n \to \infty}\int_{\T^d} \tilde K(w_n^\omega) \,\varphi \, d\t  = \xi(\varphi)$ for every continuous function $\varphi$ on $\T^d$.

This property is equivalent to (for $\PP$-a.e. $\omega$): 
\begin{eqnarray}
\lim_{n \to \infty} \int \tilde K(w_n^\omega)(\t) \, e^{-2\pi i \langle \p, \t \rangle} \ d\t 
= \lim_{n \to \infty} \int \tilde K(w_n^\omega)(\t) \, \cos (2\pi \langle \p, \t \rangle) \ d\t = \hat \xi (\p), \, \forall \p \in \Z^d. \label{coefF}
\end{eqnarray}
Another equivalent formulation is 
\begin{eqnarray}
\lim_{n \to \infty} {V_{n}(\omega, \p) \over V_{n}(\omega)} &=& \hat \xi(\p), \, \forall \p \in \Z^d, \text{ for a.e. } \omega, \nonumber\\
\text{ with } V_{n}(\omega, \p) &:=& \#\{0 \leq u, v < n: \ Z_u(\omega) - Z_v(\omega) = \p\}, \, \p \in \Z^d. \label{vnp}
\end{eqnarray}
For $f$ satisfying (\ref{sumcorr1}), it implies that the (asymptotic) normalised variance is, for a.e. $\omega$,
\begin{eqnarray}
\sigma^2(f) := \lim_n {\|\sum_{\el \in \Z^d} w_n(\el, \omega) \, T^\el f\|_2^2 \over \sum_{\el \in \Z^d} |w_n(\el,  \omega)|^2}  
= \lim_n \int_{\T^d} \tilde K(w_n^\omega)(\t) \, \varphi_f(\t) \, d\t = \xi(\varphi_f). \label{asymptVar0}
\end{eqnarray}
It can be shown that every summation associated to a random walk in $\Z^d$ is $\xi$-regular for some measure $\xi$ (cf. \cite{CohCo17}).

We summarize below the results on the asymptotic variance (see \cite{CohCo17} for the proofs).

\subsubsection{\bf Recurrence/transience}

\

Recall that a r.w. $Z = (Z_n)$ is recurrent if and only if $\sum_{n=1}^\infty \PP(Z_n = \0) = + \infty$. 

Let $m_1(Z):=\sum_{\el \in \Z^d} \, \PP(\zeta_0 = \el) \, \|\el\|$,  $m_2(Z):=\sum_{\el \in \Z^d} \, \PP(\zeta_0 = \el) \, \|\el\|^2$.
 
For $d=1$, if $m_1(Z)<\infty$, then $Z$ is recurrent if and only if it is centered;
\hfill \break for $d=2$, if $m_2(Z)<\infty$, then $Z$ is recurrent if and only if it is centered;
\hfill \break for $d\ge3$, if $m_2(Z)<\infty$, then it is always transient.

We denote by $\Psi(\t) = \E [e^{2\pi i \langle \zeta_0, \t\rangle}], \ \t \in \T^d$, the characteristic function of the r.w. and put
\begin{eqnarray}
\Phi(\t) := {1 - |\Psi(\t)|^2 \over |1 - \Psi(\t)|^2}. \label{defWt}
\end{eqnarray}
\begin{rem} \label{nonnull} For $\t \not = \0$ in $\T^d$, $\Phi$ is well defined (since $Z$ is aperiodic), nonnegative and $\Phi(\t) = 0$ only on $\Gamma_1 \setminus \{\0\}$. 
Hence it is positive for a.e. $\t$, except when the r.w. is `deterministic'
(i.e., if $\PP(\zeta_0 = \el) = 1$ for some $\el \in \Z^d$, so that $|\Psi(\t)| \equiv 1$ in this case).
\end{rem}
A r.w. of genuine dimension $d$ which is aperiodic is transient or recurrent depending on whether $\Re e ({1 \over 1- \Psi})$ is integrable or not 
on the $d$-dimensional unit cube (\cite{Sp64}).

\goodbreak
{\bf Transient case}

In the transient case, one can show:
\begin{thm} \label{transrecSpec} (\cite{Sp64}) Let $Z =(Z_n)$ be a transient aperiodic random walk in $\Z^d$.

a) The function $\Phi$ is integrable on $\T^d$ and, with a nonnegative constant $K$, we have
\begin{eqnarray*}
&&I(\el) := 1_{\el= \0} + \sum_{k = 1}^{\infty} \, [\PP(Z_k = \el) + \PP(Z_k = - \el)] = \int_{\T^d} \cos (2\pi \langle \el, \t \rangle) \, \Phi(\t) \, d\t + K,
\forall \el \in \Z^d.
\end{eqnarray*}
$b_1$) Suppose $d = 1$. If $m_1(Z) = +\infty$, then $K =0$. 
If $m_1(Z) < \infty$, then $Z$ is non centered (because it is transient) and $K = |\sum_{\ell \in \Z} \, \PP(X_0 = \ell) \, \ell|^{-1}$. 

$b_2$) If $d > 1$, then $K = 0$.

c) Denoting by $d\xi(\t)$ the measure $\Phi(\t) d\t + K \delta_\0(\t)$, we have, for a.e. $\omega$, 
\begin{eqnarray}
&&\int \frac1n K_n^\omega(\t) \cos (2\pi \langle \el, \t\rangle) \, d\t =  1_{\el= \0} +
\frac1n \, \sum_{k = 1}^{n-1} \, \sum_{j = 0}^{n-k-1} [1_{Z_k (\theta^j \omega) = \el} + 1_{Z_k (\theta^j \omega) = - \el}] \nonumber \\
&& \underset{n \to \infty} \to I(\el) = \int  \cos (2\pi \langle \el, \t\rangle) \,  d\xi(\t). \label{convFour}
\end{eqnarray}
\end{thm} 
It follows that the summation along a transient r.w. behaves for the normalisation like the iteration of a single transformation, is $\xi$-regular (up to a constant factor) and that
\begin{eqnarray}
\lim_n \frac1n \, \|\sum_{k=0}^{n-1} T^{Z_k(\omega)} f\|_2^2 = \int \Phi(\t) \, \varphi_f(\t) \, d\t + K \varphi_f(\0). \label{VarTrans1}
\end{eqnarray}
From (\ref{convFour}) , (\ref{VarTrans1}) and the expression of $\varphi_f$, $\varphi_f(\t) = \sum_{\el \in \Z^d} \langle T^\el f, f\rangle \cos (2\pi \langle \el, \t\rangle)$, 
we deduce:
\begin{eqnarray*}
&&\int \Phi(\t) \, \varphi_f(\t) \, d\t + K \varphi_f(\0) = \int \varphi_f(t) \, d\xi(\t) = 
\sum_{\el \in \Z^d} \langle T^\el f, f\rangle \int \cos (2\pi \langle \el, \t\rangle) d\xi(\t)= \\
&&\sum_{\el \in \Z^d} I(\el) \, \langle T^\el f, f\rangle =
\|f\|_2^2 + 2 \sum_{k \geq 1} \bigl(\sum_{\el \in \Z^d}  \PP(Z_k = \el) \, \langle T^\el f, f\rangle \bigr).
\end{eqnarray*}
\begin{rem} (about the variance in the non deterministic transient case)

Let $f$ be in $L^2(E, \mu)$ with real values and satisfying (\ref{sumcorr1}). By Remark \ref{nonnull}, the (quenched) asymptotic variance is $\not = 0$, if $f$ is not a.e. equal to 0. 
Let $F_f(\omega, x) = F(\omega, x) = f(x)$.
For the map $T_{\zeta_0}$ acting on the product space $\Omega \times E$ endowed with the product measure $\PP \times \mu$, we have 
$$\int T_{\zeta_0}^n F \, F d\mu \,d\PP = \sum_\el (\E 1_{Z_n = \el}) \, \langle T^\el f \, f \rangle = \sum_{\el \in \Z^d}  \PP(Z_n = \el) \, \langle T^\el f, f\rangle, \, n \geq 0.$$
In the transient case it holds, for every $\el$, $\sum_{k \geq 1} \PP(Z_k = \el) \leq \sum_{k \geq 1} \PP(Z_k = \underline 0) < +\infty$.
Therefore the density of the spectral measure for $F_f$ and the map $T_{\zeta_0}$ is
$$\|f\|_2^2 + 2 \sum_{k \geq 1} \bigl(\sum_{\el \in \Z^d}  \PP(Z_k = \el) \, \langle T^\el f, f\rangle \bigr) \, \cos 2\pi k t.$$
The asymptotic variance $\lim_n \frac1n \, \|\sum_{k=0}^{n-1} T_{\zeta_0}^k F_f\|_2^2$ for the annealed model is the same as for the quenched model and is equal to
$\displaystyle \|f\|_2^2 + 2 \sum_{k \geq 1} \bigl(\sum_{\el \in \Z^d}  \PP(Z_k = \el) \, \langle T^\el f, f\rangle \bigr)$.

It follows that the function $F_f$ on $\Omega \times E$ (which depends only on the second coordinate), 
with $f$ as above and non a.e. null, is never a coboundary in $L^2(\PP \times \mu)$ for $T_{\zeta_0}$, 
because the asymptotic variance is non null. Observe also that $F_f$ even is not a measurable coboundary, at least when the CLT holds, 
which is the case of the situations that we are going to considered here. This follows from the fact that, for a single measure preserving transformation, 
if an observable is a coboundary in the space of measurable functions, then the limiting distribution of the ergodic sums after normalisation 
by any sequence tending to infinity is the Dirac mass at 0, which is excluded here.
\end{rem}

\vskip 3mm
{\bf Recurrent case}

Let us consider now the case $d=2$ and a centered random walk $Z$ with a moment of order 2.
By the local limit theorem (LLT), $Z$ is recurrent.

A non standard normalization occurs in the CLT for sums along $Z_n$ as recalled below.
There are $C_0, C$ finite positive constants
\footnote{If the r.w. is strongly aperiodic, $C_0 = (\pi \sqrt{\det \Sigma})^{-1}$.}
such that (cf. \cite{Bo89} Lemma 2.6, \cite{Lew93} Proposition 1.4 for (\ref{BoundVar}) and (\ref{equivAn20}), \cite{CohCo17} Theorem 4.13 for (\ref{equivAn2})):
\begin{eqnarray}
&&\E(V_n) \sim C_0 n \ln n, \ \Var(V_n) \leq C n^2, \label{BoundVar}\\
&&\varphi_n(\omega) := {V_n(\omega) \over C_0 n \ln n} \to 1, \text{ for a.e. } \omega, \label{equivAn20} \\ 
&&\varphi_n(\omega, \p) := {V_n(\omega, \p) \over C_0 n \ln n}  \to 1, \forall \p \in\Z^d, \text{ for a.e. } \omega. \label{equivAn2}
\end{eqnarray}
Therefore the summation along the r.w. $Z$ is $\delta_\0$-regular: the normalised kernel satisfies 
$\lim_n \int \tilde K(w_n^\omega)(\t) \, e^{-2\pi i \langle \p, \t \rangle} \ d\t = 1, \, \forall \p \in \Z^d$ and the asymptotic variance is
\begin{eqnarray}
\sigma^2(f) = && \lim_n \, (C_0 n \ln n)^{-1} \|\sum_{k=0}^{n} T^{Z_k(\omega)} f\|_2^2 = \sum_{\k \in \Z^d} \langle T^\k f \, f \rangle = \varphi_f(\0). \label{omegaReg1}
\end{eqnarray}

The results presented below are valid for the cases covered above, hence excludes only the one-dimensional recurrent case.

We stress that, in the recurrent 2-dimensional case, the variance can be degenerate, while this does not occur in the transient case unless $f = 0$.

\subsubsection{\bf Number of self-intersections of a 2-dimensional centered r.w.} \label{selfInter}

\

In this subsection, we study more precisely the case $d=2$ and a centered random walk $(Z_n)$ with a moment of order 2.  

If $I, J$ are intervals, the quantity $V(\omega, I, J, \p) :=$
\begin{eqnarray}
&&\int \bigl(\sum_{u \in I}  e^{2\pi i \langle Z_u(\omega), \t\rangle}\bigr) 
\, \bigl(\sum_{v \in J} \, e^{- 2\pi i \langle Z_v(\omega), \t \rangle}\bigr) \, e^{-2\pi i \langle \p, \t\rangle} \, d\t 
\ = \#\{(u, v) \in I \times J: \, Z_u(\omega) - Z_v(\omega) = \p\} \label{cross1}
\end{eqnarray}
is non negative and increases when $I$ or $J$ increases for the inclusion order.

We write simply $V(\omega, I, \p)$ if $I=J$,  $V(\omega, I)$ for  $V(\omega, I, \0)$, $ V_n(\omega)$ and $V_n(\omega, \p)$ as above for $V(\omega, [0, n[)$
and  $V(\omega, [0, n[, \p)$.

Observe that $V(\omega, J) = \sum_{\el \in \Z^2} w(\omega, J, \el)^2$, where $w(\omega, J, \el) = \sum_{i \in J} 1_{Z_i(\omega) = \el}$. 
Notice also that $V(\omega, [b, b+k[) = V(\theta^b\omega, [0, k[) = V_k(\theta^b\omega)$, for $b \geq 0, k \geq 1$.

Let $A, B$ be in $[0, 1]$. For simplicity, in the formulas above and below, we write $nA$, $nB$ instead of $\lfloor nA \rfloor$ or $\lfloor nA \rfloor +1$, 
$\lfloor nB \rfloor$, $\theta^t$ instead of $\theta^{\lfloor t \rfloor}$. The equalities are satisfied up to the addition of quantities which are bounded independently from $A, B, n$.
We have:
\begin{flalign*}
&V(\omega, [nA, nB], \p) =  \int (\sum_{u \in [nA, nB]}  e^{2\pi i \langle Z_u(\omega), \t\rangle}) 
\, (\sum_{v\in [nA, nB]} \, e^{- 2\pi i \langle Z_v(\omega), \t \rangle}) \, e^{-2\pi i \langle \p, \t\rangle} \, d\t& \\
&= \#\{u, v \in [0, n(B-A)]:\, \sum_{i=0}^{u-1} \zeta_0(\theta^{i+nA} \omega) - \sum_{i=0}^{v-1} \zeta_0(\theta^{i+nA} \omega) = \p\}
= V(\theta^{nA}\omega, [0, n(B-A)], \p).&
\end{flalign*}
By (\ref{equivAn20}) and (\ref{equivAn2}) there a set $\Omega_0$ of full probability such that
\begin{eqnarray}
&&V_n(\omega) \leq K(\omega) \, n \ln n, \, \forall n \geq 2, \text{ where the function } K \geq 0 \text{ is finite on } \Omega_0, \label{Komega}\\
&& \text{ for any fixed } A \in ]0, 1], \ V(\omega, [1, nA], \p) \sim C_0 n A \ln n, \text{ for } \omega \in \Omega_0. \label{equivAn}
\end{eqnarray}
By \cite[Lemma 2.5]{Bo89} we have
\begin{eqnarray}
\sup_{\el\in \Z^2} \, w_n(\omega, \el)=o(n^\varepsilon), \text{ for a.e. } \omega , \text{ for every } \varepsilon>0. \label{nVarep1}
\end{eqnarray}
For a simple r.w. on $\Z^2$, Erd\"os and Taylor (\cite{ErTa60}) have shown:
$\displaystyle \limsup_n \, \sup_{\el\in \Z^2} \, {w_n(\omega, \el) \over (\log n) ^2} \leq \frac1\pi$.

The result has been extended by Dembo, Peres, Rosen and Zeitouni who proved for an aperiodic centered random walk on $\Z^2$ with moments of all orders \cite{DPRZ01}:
\begin{eqnarray} 
\lim_n \sup_{\el\in \Z^2} {w_n(\omega, \el) \over (\log n) ^2} = \frac1\pi. \label{DPRZbound}
\end{eqnarray}

We will need also to bound, for $\el_1, \el_2, \el_3 \in \Z^d$, $W_n(\omega, \el_1, \el_2, \el_3) :=$
\begin{equation} 
\#\{1 \leq i_0, \, i_1, \, i_2, \, i_3 < n: \, Z_{i_1}(\omega) - Z_{i_0}(\omega) = \el_1, 
Z_{i_2}(\omega) - Z_{i_0}(\omega)= \el_2, Z_{i_3}(\omega) - Z_{i_0}(\omega) = \el_3\}. \label{wn123}
\end{equation}
\begin{lem} \label{majwnm} 
There exists a positive integrable function $C_3$ such that 
\begin{eqnarray} 
W_n(\omega, \el_1, \el_2, \el_3) \leq C_3(\omega) \, n \, (\ln n)^5, \, \forall n \geq 1. \label{majnumb3}
\end{eqnarray}
\end{lem}
\proof It suffices to bound the sum with strict inequality between indices
$$W_n'(\omega) = \sum_{1 \leq i_0 < i_1 < i_2 < i_3 \leq n} 1_{Z_{i_1} - Z_{i_0} = \, \el_1} . 1_{Z_{i_2} - Z_{i_1}= \, \el_2 - \el_1}. 
1_{Z_{i_3} - Z_{i_2} = \, \el_3 - \el_2}.$$
Using independence and the local limit theorem for the random walk, we find the bound
\begin{eqnarray*} 
\int W_n'(\omega) \, d\PP(\omega)  \leq  C_1 \sum_{i_0, i_1, i_2, i_3 \in [1, n]} (i_1 \, i_2 \, i_3)^{-1} \leq C_2 \, n \, (\ln n)^3.
\end{eqnarray*}
Therefore $\displaystyle \sum_{p=1}^\infty \int 2^{-p}(\ln(2^p))^{-5} \, W_{2^p}' \, d\PP < \infty$.
The function $C(\omega) := \sum_{p=1}^\infty 2^{-p}(\ln(2^p))^{-5} \, W_{2^p}'$ is integrable and we have:
$W_{2^p}'(\omega) \leq C(\omega) \, 2^p \, (\ln(2^p))^5, \, \forall p \geq 1$.

Let $p_n$ be such that: $2^{p_n-1} \leq n < 2^{p_n}$. Since $W_n'$ is increasing with $n$, we obtain:
$$W_{n}'(\omega) \leq W_{2^{p_n}}'(\omega) \leq C(\omega) \, 2^{p_n} \, (\ln (2^{p_n}))^5 \leq C(\omega) \, 2n \, (\ln (2n))^{5}
\leq C'(\omega) \, n (\ln n)^{5}. \eop$$

{\bf Variance for the finite dimensional distributions}

The following lemma will be applied to the successive return times of a point $\omega$ into a set under the iteration of the shift $\theta$.
\begin{lem} \label{mean0} Let $(y(j), j \geq 1)$ be a sequence with values in $\{0, 1\}$ such that $\lim_n \frac1n \sum_{j=1}^n y(j) = a > 0$.
If $(k_r)$ is the sequence of successive times such that $y(k_r) = 1$, then, for every $\delta >0$, there is $n(\delta)$ such that,
 for $n \geq n(\delta)$, $k_{r+1} - k_r \leq \delta n$, for all $r \in [1, n]$.
\end{lem}
\Proof Since $r = \sum_{j=1}^{k_r} \, y(j)$, we have: $k_r / r = k_r / \sum_{j=1}^{k_r} \, y(j)  \to a^{-1}$. 
Hence, for every $\delta > 0$, there is $n_1(\delta)$ such that $0 < k_{r+1}  - k_r \leq \delta r$, for $r \geq n_1(\delta)$.
Therefore, if $n \geq n_1(\delta)$, then  $0 < k_{r+1}  - k_r \leq \delta r \leq \delta n$, for $r \in [n_1(\delta), n]$. 

If $n(\delta) \geq n_1(\delta)$ is such that $k_{r+1}  - k_r \leq \delta n(\delta)$ for $r \leq n_1(\delta)$, we get the result of the lemma. \eop

\begin{lem} \label{mean1} Let $\Lambda$ be a measurable set in $\Omega$ of positive measure. 
Let $k_r= k_r(\omega)$ be the successive times such that $\theta^{k_r} \omega \in \Lambda$.
For a.e. $\omega$, for every positive small enough $\delta$, there is $n(\delta)$ such that for $n \geq n(\delta)$
\hfill \break
1) $k_{r+1} - k_r \leq \delta n$, for all $r \in [1, n]$; moreover, $k_n \sim c n$, with $c=  \PP(\Lambda)^{-1}$, when $n \to \infty$;
\hfill \break
2) there are integers $v < 2/\delta$ and $0 = \rho_1^{(n)} < \rho_2^{(n)} < ... < \rho_v^{(n)} \leq n < \rho_{v+1}^{(n)}$, 
such that $\theta^{\rho_i^{(n)}} \omega \in \Lambda$ and $\frac12 \delta n \leq \rho_{i+1}^{(n)} - \rho_i^{(n)} \leq \frac32 \delta n$, for $i = 1, ..., v$.
\end{lem}
\Proof Since $\theta$ is ergodic on $(\Omega, \PP)$, Birkhoff ergodic theorem implies
$\lim_n \frac1n \sum_0^{n-1} 1_\Lambda (\theta^k \omega) = \PP(\Lambda) > 0$, for a.e. $\omega$ and $k_n / n \to \PP(\Lambda)^{-1}$.
Hence Lemma \ref{mean0} implies 1). For 2), we select in the sequence $(k_r)$ an increasing sequence of visit times to the set $\Lambda$ satisfying the prescribed conditions 
by eliminating successive times which are at a distance $< \frac12 \delta n$. \eop

\vskip 2mm
{\it Asymptotic orthogonality of the cross terms}
 
\begin{proposition} \label{asympOrth} For $0 < A < B < C < D < 1$, $\p \in \Z$,
\begin{eqnarray}
&&\int \bigl[(\sum_{v = nA}^{nB}  e^{2\pi i \langle Z_v(\omega), \u \rangle}) 
\, (\sum_{w = nC}^{nD} \, e^{- 2\pi i \langle Z_w(\omega), \u \rangle}) \, e^{-2\pi i \langle \p, \u \rangle} \bigr] \, d\u
=\varepsilon_n(\omega) \, n \log n, \text{ with } \varepsilon_n(\omega) \to 0. \label{orth1}
\end{eqnarray}
\end{proposition}

The above integral is the non negative self-intersection quantity: $V(\omega, [nA, nB], [nC, nD], \p)$. By (\ref{cross1}), $V(\omega, I, J, \p)$ increases when $I$ or $J$ increases. 
Hence, it suffices to show (\ref{orth1}) for the intervals $[1, nA], [nA, n]$, for $0 < A < 1$. The proof below is based on (\ref{equivAn2}) and (\ref{equivAn}).

\begin{lem} \label{asympOrtho1} There is a set $\hat \Omega \subset \Omega$ such that $\PP(\hat \Omega) = 1$ and for all $\omega \in \hat \Omega$,
the following holds:
\begin{eqnarray}
&&\lim_n \varphi_{n B}(\theta^{n A}\omega, \, \p) = \lim_n {V(\omega, [nA, n], \p) \over C_0 \, n B \ln n} = 1, \text{ for } A\in ]0,1[, B = 1 - A; \label{lowup0} \\
&&V(\omega, [1, nA], [nA, n], \p) + V(\omega,[nA, n], [1, nA], \p) =\varepsilon_n(\omega) \, n \log n, \text{ with } \varepsilon_n(\omega) \to 0. \label{orth2}
\end{eqnarray}
\end{lem}
\proof  {\it 1) The set $\hat \Omega$.} 
For every $L \geq 1$ and $\delta > 0$, let $\Lambda(L, \delta):= \{\omega: \varphi_n(\omega, \, \p)  - 1 \in [-\delta, \delta], \forall n \geq L\}$.
We have $\lim_{L \uparrow \infty} \PP(\Lambda(L, \delta)) = 1$. There is $L(\delta)$ such that $\PP(\Lambda(L(\delta), \delta)) \geq \frac12$.

Let $(\delta_j)$ be a sequence tending to 0. We apply Lemma \ref{mean1} to $\Lambda(L(\delta_j), \delta_j)$ for each $j$.
By taking the intersection of the corresponding sets, we get a set of $\omega$'s of full $\PP$-measure. The set $\hat \Omega$ is the intersection of this set 
with the set $\Omega_0$ (of full measure) for which the law of large numbers holds for $(V_n(\omega))$. Let $\omega \in \hat \Omega$. 

{\it 2) Proof of (\ref{lowup0})}. 
We have $V(\omega, [nA, n[, \p) = V(\theta^{nA} \omega, [0, n(1-A)[, \p)$ and
\begin{eqnarray}
&&V(\omega, [1, n], \p) - V(\omega, [1, nA[, \p) - V(\omega, [nA, n], \p) \nonumber \\
&& = V(\omega, [1, nA[, [nA, n[, \p) + V(\omega,[nA, n], [1, nA[, \p) \geq 0. \label{differ1}
\end{eqnarray}
Claim: for an absolute constant $C_1$ depending on $A$ and $\p$, for every $\delta$, for $n$ big enough,
\begin{eqnarray}
\varphi_{n B}(\theta^{n A}\omega, \, \p) = {V(\omega, [nA, n], \p) \over C_0 \, n (1 - A) \ln n} \in [1 - C_1 \delta, 1 + C_1 \delta]. \label{lowup1}
\end{eqnarray}

{\it Upper bound:} The law of large numbers for $V_n(\omega, \p)$ implies, with $|\varepsilon_n|, |\varepsilon_n'| \leq \delta$ for $n$ big enough,
$$C_0^{-1} \, V(\omega, [1, n], \p) =  (1 + \varepsilon_n) \, n \ln n, \ C_0^{-1} \, V(\omega, [1, nA], \p) = (1 + \varepsilon_n') \, nA \ln n.$$
With $B= 1 - A$, this implies by (\ref{differ1}) 
\begin{eqnarray*}
{V(\omega,  [nA, n], \p) \over C_0 \, nB \ln n} \leq {(1+ \varepsilon_n) \, n \ln n - (1+ \varepsilon_n') \, n A \ln n \over nB \ln n} 
\leq 1 + {|\varepsilon_n| \over B} + {|\varepsilon_n'| A \over B} \leq 1+ {1 + A \over B} \delta.
\end{eqnarray*}

{\it Lower bound:} We apply Lemma \ref{mean1} to $\Lambda(L(\delta), \delta)$.
 Let $n_A, n_A'$ be two consecutive visit times $\leq n$ such that $n_A \leq n A < n_A'$. 
For $n$ big enough, we have $0 <  n_A' -  n_A \leq \delta n$ and 
$$n_A = n A \, (1 - \rho_n), \ n_A' = n A \, (1 + \rho_n'), \text { with } 0 \leq A \rho_n, A \rho_n'  \leq \delta.$$
Since $\omega \in \hat \Omega$, we have for $n$ big enough, with $|\delta_n'| \leq \delta$,
\begin{eqnarray*}
&C_0^{-1} \, V(\omega, [n'_{A}, n], \p) \geq (1-\delta_n') (n - n_A') \,\ln(n- n_A') = (1-\delta_n')  (nB -nA \rho_n') \,\ln(nB -nA \rho_n') .
\end{eqnarray*}
It follows, for $\delta$ (hence $\rho_n'$) small:
\begin{flalign*}
&{V(\omega, [n_A', n], \p) \over C_0\,  (1-\delta_n')  \, n B \ln (n B)} \geq {(nB -nA \rho_n') \,\ln(nB -nA \rho_n') \over n B \ln (n B)} 
= {(B -A \rho_n') \, [\ln(nB) + \ln (1 - \frac AB \rho_n')] \over B \ln (n B)}& \\
&\geq (1- \frac AB  \rho_n') - 2 (1- \frac AB \rho_n')\frac{\frac AB \rho_n'}{\ln(nB)} 
\geq 1- \frac AB \rho_n'- 2 \frac{\frac AB \rho_n'}{\ln(nB)} \geq 1-B^{-1}\delta (1 + \frac{2}{\ln(nB)}).&
\end{flalign*}
As $V(\omega, J, \, \p)$ increases when the set $J$ increases, we have by the choice of $n_A$ and $n_A'$:
$$V(\omega, [n_A', n], \p) \leq V(\omega, [nA, n], \p).$$ 
Therefore, for $n$ such that $\ln (nB) \geq 2$, we have
$${V(\omega, [nA, n], \p) \over C_0 \, n B \ln (n B)} \geq (1-\delta) \, (1- {2 \over B}\delta) \geq 1 - \delta (1 + {2 \over B}).$$
This shows the lower bound. Altogether with the upper bound, this proves the claim (\ref{lowup1}).

{\it 3)  Proof of (\ref{orth2}).} Let $\delta > 0$. According to (\ref{differ1}) and (\ref{lowup1}), for $n$ big enough, we have with $|\varepsilon_n''| \leq C_1 \delta$:
\begin{flalign*}
&V(\omega, [1, nA], [nA, n], \p) + V(\omega,[nA, n], [1, nA], \p) = V(\omega, [1, n], \p) - V(\omega, [1, nA], \p) - V(\omega, [nA, n], \p) &\\
&=  C_0[ (1 + \varepsilon_n) \, n \ln n  - (1 + \varepsilon_n') \, n A \ln n - (1 + \varepsilon_n'') \, n(1-A) \ln n  \leq (2 + C_1) \, C_0 \, \delta \, n \ln n.  \eop&
\end{flalign*}
Let $a_1, ..., a_s$ be real numbers and $0 = t_0 < t_1 < ... < t_{s-1} < t_s = 1$ a subdivision of $[0, 1]$.
For the asymptotic variance of $\sum_{j= 0}^s \, a_j \sum_{k=nt_{j-1} }^{nt_j} T^{Z_k(\omega)} f$, which is used later, we need the following lemma.  
Recall that $f$ has a continuous spectral density $\varphi_f$.
\begin{lem} \label{VarrwLem1} For a.e. $\omega$ and for every partition $(t_j)$, we have
\begin{eqnarray}
(C_0 n \ln n)^{-1} \|\sum_{j= 1}^s \, a_j \sum_{k=nt_{j-1}}^{nt_j} T^{Z_k(\omega)} f\|_2^2
\to \varphi_f(\underline 0)\sum_{j=1}^s a_j^2(t_j-t_{j-1}). \label{VarMulti1}
\end{eqnarray}
\end{lem}
\proof 1) Recall that proving (\ref{VarMulti1}) amounts to prove
$$(C_0 n \ln n)^{-1}  \int |\sum_{j= 1}^s \, a_j \sum_{k=nt_{j-1}}^{nt_j} e^{2\pi i \langle Z_k(\omega) , \u \rangle}|^2 \, \varphi_f(\u) \, d\u
\to \varphi_f(\underline 0)\sum_{j=1}^s a_j^2(t_j-t_{j-1}).$$

1) First suppose that $\varphi_f$ is a trigonometric polynomial $\rho$, which allows to use (\ref{orth1}) for a finite set of characters $e^{-2\pi i \langle \p, \u \rangle}$.
Using (\ref{omegaReg1}) for the asymptotic variance starting from 0, we have 
$(C_0 n \ln n)^{-1} \|\sum_{k=0}^{\lfloor tn \rfloor} T^{Z_k(\omega)} f\|_2^2 \to  t \rho(\0)$, for $t \in ]0, 1[$.
By Lemma \ref{asympOrtho1}, 
\begin{eqnarray*}
&&(C_0 n \ln n)^{-1} \|\sum_{k=\lfloor s n \rfloor}^{\lfloor t n \rfloor} T^{Z_k(\omega)} f\|_2^2 \to  (t-s) \, \rho(\0), \text{ for } 0 < s < t < 1.
\end{eqnarray*}
Expanding the square and using that the cross terms are asymptotically negligible, we have
\begin{eqnarray*}
&&(C_0 n \ln n)^{-1} \int |\sum_{j= 1}^s \, a_j \sum_{k=nt_{j-1}}^{nt_j} e^{2\pi i \langle Z_k(\omega) , \u \rangle}|^2 \, \rho(\u) \, d\u \\
&& \sim \ (C_0 n \ln n)^{-1} \bigl(\sum_{j=1}^s a_j^2 \,  \int |\sum_{k=nt_{j-1}}^{nt_j} e^{2\pi i \langle Z_k(\omega) , \u \rangle}|^2 \, \rho(\u) \, d\u \bigr)
\to \rho(\0) \sum_{j=1}^s a_j^2(t_j-t_{j-1}).
\end{eqnarray*}
This shows (\ref{VarMulti1}) for trigonometric polynomials. 

2) For a general continuous spectral density $\varphi_f$, for $\varepsilon >0$, 
let $\rho$ be a trigonometric polynomial, such that  $\|\varphi_f - \rho\|_\infty < \varepsilon$. Remark that
$$\int |\sum_{j= 1}^s \, a_j \sum_{k=nt_{j-1}}^{nt_j} e^{2\pi i \langle Z_k(\omega) , \u \rangle}|^2 \, d\u 
\leq \sum_{j, j'= 1}^s \,a_j a_{j'} \, V(\omega, [nt_{j-1}, nt_j],  [nt_{j'-1}, nt_{j'}], \0) \leq (\sum_{j= 1}^s \,|a_j|)^2 \, V_n(\omega).$$
Therefore we have:
\begin{eqnarray*}
&&\bigl|(C_0 n \ln n)^{-1} \int |\sum_{j= 1}^s \, a_j \sum_{k=nt_{j-1}}^{nt_j} e^{2\pi i \langle Z_k(\omega) , \u \rangle}|^2 \, \varphi_f(\u) \, d\u
- \varphi_f(\0) \sum_{j=1}^s a_j^2(t_j-t_{j-1})\bigr|\\
&&\leq \bigl|(C_0 n \ln n)^{-1} \int |\sum_{j= 1}^s \, a_j \sum_{k=nt_{j-1}}^{nt_j} e^{2\pi i \langle Z_k(\omega) , \u \rangle}|^2 \, \rho(\u) \, d\u
- \rho(\0) \sum_{j=1}^s a_j^2(t_j-t_{j-1})\bigr| \\
&&+ \varepsilon \, [(C_0 n \ln n)^{-1} \int |\sum_{j= 1}^s \, a_j \sum_{k=nt_{j-1}}^{nt_j} e^{2\pi i \langle Z_k(\omega) , \u \rangle}|^2 \, d\u
+ \sum_{j=1}^s a_j^2(t_j-t_{j-1})].
\end{eqnarray*}
By the remark, the above quantity inside $[\ ]$ is less than $(\sum_{j= 1}^s \,|a_j|)^2 \, (C_0 n \ln n)^{-1} V_n(\omega) + \sum_{j=1}^s a_j^2(t_j-t_{j-1})$,
which is bounded uniformly with respect to $n$. Therefore we can conclude for a general continuous spectral density by step 1).\eop

\begin{rems} \label{general_0} 1)  In Lemma \ref{mean1}, the dynamical system $(\Omega, \theta, \PP)$ can be replaced by any ergodic dynamical system.

2) If the spectral density is constant (i.e., when the $X_\k$'s are pairwise orthogonal), (\ref{orth1}) and (\ref{VarMulti1}) are a consequence of 
the law of large numbers for the number of self-intersections, that is ${V_n(\omega) \over C_0 n \ln n} \to 1$. 
The law of large numbers for $V_n(\omega, \p)$, $\p \not = \0$, is not needed. 
\end{rems}

3) A result analogous to Proposition \ref{asympOrth} is valid if the r.w. $Z$ is transient: for $0 < A < B < C < D < 1$, $\p \in \Z$,
\begin{eqnarray}
&&\int \bigl[(\sum_{v = nA}^{nB}  e^{2\pi i \langle Z_v(\omega), \u \rangle}) 
\, (\sum_{w = nC}^{nD} \, e^{- 2\pi i \langle Z_w(\omega), \u \rangle}) \, e^{-2\pi i \langle \p, \u \rangle} \bigr] \, d\u
=\varepsilon_n(\omega) \, n, \text{ with } \varepsilon_n(\omega) \to 0. \label{orthtrans}
\end{eqnarray}

4) The quenched FCLT shown in the different examples below is valid for a set of $\omega$'s of $\PP$-measure 1 given by the results of  
this Section \ref{varGene1}. This set does not depend on the $\Z^d$-dynamical systems considered in the further sections. 
The joint distribution on $\Omega \times E$ is used only when the annealed model is mentionned, like for Bolthausen's result recalled in Section \ref{indepSec}.

\subsubsection{\bf Formulation of the quenched FCLT}

\

Let $(Y_n(t), t \in [0, 1])$ be a process on $(E, \mu)$ with values in the space  $C[0, 1)$ of real valued continuous functions on $[0,1]$ or in the space $D[0, 1]$ 
of right continuous real valued functions on $[0,1]$ with left limits, endowed with the uniform norm. 

Let $(W(t), t \in [0, 1])$ be the Wiener process on $[0, 1]$. To show a functional limit theorem (FCLT) for $(Y_n(t), t \in [0, 1])$, 
i.e., weak convergence to the Wiener process, it suffices to prove the two following properties (where ``$\Longrightarrow$'' denotes the convergence in distribution):

{\it 1) Convergence of the finite dimensional distributions:} 
$$\forall \, 0=t_0 < t_1 <... < t_r =1, \ (Y_{n}(t_1), ..., Y_{n}(t_r)) \underset {n \to \infty} \Longrightarrow (W_{t_1}, ..., W_{t_r}),$$ 
a property which follows (by the Cram\'er-Wold device) from  
\begin{eqnarray}
&\sum_{j=1}^r a_j (Y_n(t_j) - Y_n(t_{j-1})) \Longrightarrow {\Cal N}(0, \sum_{j = 1}^r a_j^2 (t_j - t_{j-1})), \, \forall (a_j)_ {1 \leq j \leq r} \in \R. \label{Cramer0}
\end{eqnarray}

{\it 2) Tightness of the process.} The condition of tightness reads:
\begin{eqnarray}
&&\forall \varepsilon > 0, \ \lim_{ \delta \to 0^+} \limsup_n \mu(x \in E: \, \sup_{|t'-t| \leq \delta}|Y_n(\omega, x, t') - Y_n(\omega, x, t)| \geq \varepsilon) = 0. \label{tight0}
\end{eqnarray}
Let $(Z_n)$ be a random walk on $(\Omega, \PP)$ with values in $\Z^d$, $d \geq 1$, and let $X = (X_\el(x))_{\el \in \Z^d} = (T^\el f(x))_{\el \in \Z^d}$ 
be a $d$-dimensional random field defined on a probability space $(E, \cal B, \mu)$. 
A quenched FCLT is satisfied by the sums along $Z_n$ if, for $\PP$-a.e. $\omega$, the functional central limit theorem holds for the process (cf. Notation (\ref{defsumRW}))
\begin{eqnarray}
 (Y_n(\omega, x, t))_{t \in [0, 1]} := \bigl({S_{[n t]}^{\omega, X}(x) \over \sqrt{C_0 n \log n}}\bigr)_{t \in [0, 1]}. \label{defYn0}
\end{eqnarray}

\section{\bf Random walk in random scenery} \label{indepSec}

We consider in this section $d= 2$ and the random walk in random scenery $S_{n}^{\omega, X}(x)$, 
that is the process $(Y_n)$ defined by (\ref{defYn0}) when $X$ is a 2-dimensional random field of i.i.d. real variables with $\E(X_0^2) = 1$ 
and mean 0 on a probability space $(E, \cal B, \mu)$. 

It was shown by E. Bolthausen \cite{Bo89} that this process satisfies an annealed FCLT: 
with respect to the probability $\PP \times \mu$, the law of $Y_n$ converges weakly to the Wiener measure. 

We show a quenched FCLT for the r.f. $X$ (and when $X$ is a r.f. of moving averages of i.i.d. random variables in Section \ref{movAverag}).
As for the annealed FCLT in \cite{Bo89} and for Theorem 2.2 in \cite{DeUt11} for a 1-dimensional stable r.w., the proof of Theorem \ref{FCLTind1} 
below is based on the maximal inequality shown by Newman and Wright \cite{NeWr81} for associated r.v.s.
\begin{definition} {\rm (\cf \cite{EPW})
Recall that real random variables $X_1,\ldots, X_n$ are {\it associated} if, for every $n \geq 1$, for all non-decreasing (in each coordinate) functions $f,g: \R^n\mapsto \R$,
we have ${\rm Cov}(f(X_1,\ldots,X_n), \, g(X_1,\ldots, X_n))\ge 0$ (if the covariance exists).
Non-decreasing functions of a family of associated random variables are associated \cite{EPW}. 
Independent variables are associated.

It follows that, if ($X_\el, \el \in \Z^2$), are associated r.v.s, in particular independent, then the r.v.s ($X_{Z_k(\omega)}, k \geq 0$), 
are associated for every $\omega \in \Omega$.}
\end{definition}
\begin{thm} \label{FCLTind1} If $\E(X_\0^2) = 1$, for $\PP$-a.e. $\omega$, the process 
$\displaystyle \bigl (Y_n(\omega, x, t\bigr)_{t \in [0,1]} = \bigl({S_{\lfloor nt \rfloor}^{\omega, X} (x) \over\sqrt{ n \log n}}\bigr)_{t \in [0,1]}$
satisfies a FCLT with asymptotic variance $\sigma^2 = (\pi \sqrt{\det \Sigma})^{-1}.$
\end{thm}
\proof  1) For the convergence of the finite dimensional distributions, the proof, relying on Cram\'er-Wold's theorem and Lindeberg's CLT, is as in \cite{Bo89}. 
Another proof, based on truncation and cumulants, can be given, like for the more general case of moving averages in Section \ref{movAverag}.

2) {\it Tightness of the process $(Y_n)$.} The following is shown in \cite{NeWr81}, p. 673:

Let $U_1, U_2, \ldots$ be centered associated random variables with finite second order moment.
Put $S_k=\sum_{j=1}^k U_j,$ for $k\geq1$. Then, for every $\lambda>0$ and $n \geq 1$, we have
\begin{eqnarray}
\mu(\max_{1\le k \le n}|S_k| \geq \lambda \, \|S_n\|_2) \leq 2\mu \bigl(|S_n|\ge (\lambda-\sqrt{2}) \, \|S_n\|_2 \bigr). \label{NWinequal1}
\end{eqnarray}
Inequality (\ref{NWinequal1}) can be applied for every fixed $\omega$ to $U_j =  X_{Z_j(\omega)}$ and to the sums $S_J=\sum_{j=b}^{b+k} X_{Z_j(\omega)}$ 
for any interval $J = [b , b+k] \subset [0, n]$. We also note that $\E(S_J^2) = \|X_0\|_2^2 \, V(\omega, J)$.  

a) First, let us assume that $\E(X_\0^4)  < \infty$. With $K$ given by (\ref{Komega}), we have 
\begin{eqnarray}
&&\|\sum_{i \in J} X_{Z_i(\omega)}\|_{4, \mu}^4 
= 3 \E(X_{\0}^2)^2 \, \sum_{\el_1 \not = \el_2} w(\omega, J, \el_1)^2 \, w(\omega, J, \el_2)^2 + \, \E(X_\0^4) \, \sum_{\el} w(\omega, J, \el)^4 \nonumber \\
&&\ \leq 4 \, \E(X_\0^4) \, V(\omega, J)^2 \leq 4 \, \E(X_\0^4) \, (K(\theta^b \omega))^2 \, (k \ln k)^2. \label{mom4iid1}
\end{eqnarray}
Let $C_1$ be a constant $> 0$ such that $\PP\{\omega: \, K(\omega) \leq C_1\} > 0$. 
Using Lemma \ref{mean1}, for $n$ big enough and $\delta \in ]0, 1[$, there are times 
$0 = \rho_1 < \rho_2< ... < \rho_v \leq n < \rho_{v+1}$, with $v < 2/\delta$, such that $K(\theta^{\rho_i} \omega) \leq C_1$ 
and $\frac12 \delta n \leq \rho_{i+1} - \rho_i \leq \frac32 \delta n$, for $i = 1, \ldots, v$. 

Let $t_i=\frac{\rho_i}{n}$, $\lambda=\frac{\varepsilon}{\sqrt{\delta}}$, $J_i=[\rho_{i-1}, \ldots, \rho_i[$, $m_i =\frac23 (\rho_{i+1} - \rho_i) \leq \delta n$. 
There is $C$ such that, by (\ref{Komega}) and (\ref{mom4iid1}),
\begin{eqnarray}
&&\|\sum_{j=\rho_{i-1}}^{\rho_i} X_{Z_j(\omega)}\|_2 \leq C \, \|X_0\|_2 \, (n \, \delta \, \log (n \delta))^\frac12, 
\|\sum_{j=\rho_{i-1}}^{\rho_i} X_{Z_j(\omega)}\|_4 \leq C \, \|X_0\|_4 \, (n \, \delta \, \log (n \delta))^\frac12, \forall i.
\end{eqnarray}
Using (\ref{NWinequal1}), we get, with $\sigma^{(i)} := \|\sum_{j=\rho_{i-1}}^{\rho_i} X_{Z_j(\omega)}\|_2$, $\lambda_i := \varepsilon \sqrt{n\log n} / \sigma^{(i)}$, 
by Chebyshev's inequality (for a moment of order 4):
\begin{eqnarray}
&&\mu(\sup_{t_{i-1}\le s\le t_i} | \sum_{j=[t_{i-1}n\rfloor}^{\lfloor sn\rfloor} X_{Z_j (\omega)} |\geq \varepsilon \sqrt{n\log n}) 
=\mu(\max_{\rho_{i-1} \leq k \leq \rho_i} |\sum_{j=\rho_{i-1}}^{k} X_{Z_j(\omega)}|\geq \lambda_i \, \sigma^{(i)}) \nonumber \\
&&\leq 2 \mu(|\sum_{j=\rho_{i-1}}^{\rho_i} X_{Z_j(\omega)}| \geq (\lambda_i -\sqrt{2}) \, \sigma^{(i)})
\leq 2 \mu(|\sum_{j=\rho_{i-1}}^{\rho_i} X_{Z_j(\omega)}| \geq \frac12 \lambda_i \, \sigma^{(i)}) \nonumber \\
&&\leq 2 \mu(|\sum_{j=\rho_{i-1}}^{\rho_i} X_{Z_j(\omega)}| \geq \frac12 \varepsilon \sqrt{n\log n}) 
\leq 2 \, C^4 \, \|X_0\|_4^4 \, {(n \, \delta \, \log (n \delta))^2 \over \frac1{16} \varepsilon^4 (n\log n)^2} 
\leq 32 \, C^4 \, \|X_0\|_4^4 \, {\delta^2 \over \varepsilon^4}. \label{majC4}
\end{eqnarray}
We have used that $\lambda_i$ is big if $\delta$ is small. Observe now that (cf. \cite{Bi99})
\begin{eqnarray*}
&&\mu(\sup_{|t'-t|\le\delta} |Y_n(t)-Y_n(s)| \geq 3 \varepsilon) \leq
\sum_{i=1}^v \mu(\sup_{t_{i-1}\le s\le t_i} | \sum_{j=[t_{i-1}n\rfloor}^{\lfloor sn\rfloor} X_{Z_j (\omega)} |\geq \varepsilon \sqrt{n\log n}).
\end{eqnarray*}
Hence, by (\ref{majC4}) we get
$\displaystyle \mu(\sup_{|t'-t|\le\delta} |Y_n(t)-Y_n(s)| \geq 3 \varepsilon) 
\leq 32 \, C^4 \, \|X_0\|_4^4 \, \frac{2}{\delta}\frac{\delta^2}{\varepsilon^4} = 64 \, C^4 \, \|X_0\|_4^4 \, \frac{\delta}{\varepsilon^4}$.

b) Now we use a truncation. For $L>0$, let 
\begin{eqnarray*}
&&\hat X_\k^L:=X_\k \, {\bf 1}_{\{X_\k \leq L\}}- \E(X_\k \, {\bf 1}_{\{X_\k \leq L\}}), 
\ \tilde X_\k^L:=X_\k -\hat X_\k^L = X_\k \, {\bf 1}_{\{X_\k > L\}}- \E(X_\k \, {\bf 1}_{\{X_\k > L\}}), \\
&&\hat Y_n^L(t) = \frac{1}{\sqrt{C_0n\log n}}\sum_{j=0}^{\lfloor tn \rfloor} \hat X^L_{Z_j(\omega)} \text { and }
\tilde Y^L_n(t) := Y_n(t) - \hat Y_n^L(t)= \frac{1}{\sqrt{C_0 n \log n}} \sum_{j=0}^{\lfloor tn \rfloor} \tilde X^L_{Z_j(\omega)}.
\end{eqnarray*}
Since we have still sums of associated random variables, all what we have done above (including (\ref{NWinequal1}) holds for both sums, 
except that for the unbounded part of the truncation we only have a moment of order 2. 
We use Chebyshev's inequality (for a moment of order 2) to control the unbounded truncated part:
\begin{eqnarray*}
\mu(|\sum_{j=\rho_{i-1}}^{\rho_i} \tilde X^L_{Z_j(\omega)}| \geq \frac12 \varepsilon \sqrt{n\log n}) 
\leq C^2 \, \|\tilde X^L_0\|_2^2 \, {n \, \delta \, \log (n \delta) \over \frac14 \varepsilon^2 \, n\log n} 
\leq 4 C^2 \, \|\tilde X^L_0\|_2^2 \, {\delta \over \varepsilon^2}.
\end{eqnarray*}
Hence, for $n$ and $\lambda$ big enough, the sum over $i$ is comparable for some constant $C'$ with
$$C\sum_{i=1}^v \frac{ \|\tilde X^L_0\|_2^2}{\lambda^2}\leq \frac{C'}\delta \frac{\delta \|\tilde X^L_0\|_2^2 }{\varepsilon^2} 
= C' \, \varepsilon^{-2} \, \|\tilde X^L_0\|_2^2.$$
Applying the inequality $\mu(|f+g|\ge\varepsilon) \leq \mu(|f|\geq \frac\varepsilon2)+\mu(|g| \geq \frac\varepsilon2)$
to $Y_n(t) = \hat Y_n^L(t) + \tilde Y^L_n(t)$, we obtain the bound:
\begin{eqnarray*}
&&\mu (\sup_{|t'-t| \leq \delta}|Y_n(t') - Y_n(t)| \geq 3\varepsilon) \leq 2^{10} C^4 \,\frac{L^4 \delta}{\varepsilon^4} + 4 C' \, \frac{\|\tilde X^L_0\|_2^2}{\varepsilon^2}.
\end{eqnarray*}
We need, for fixed $\varepsilon > 0$, $\lim_{\delta \to 0^+} \limsup_n \mu (\sup_{|t'-t| \leq \delta}|Y_n(t') - Y_n(t)| \geq 3\varepsilon)  = 0$.

Let $\eta > 0$. First we take $L$ such that  $4 C' \, \frac{\|\tilde X^L_0\|_2^2}{\varepsilon^2} < \frac12 \eta$, then $\delta$ such that
$2^{10} C^4 \,\frac{L^4 \delta}{\varepsilon^4} \leq \frac12 \eta$. \eop

\vskip 3mm
\goodbreak
{\bf A model based on the Lorentz process} \label{LorentzFCLT}

We sketch briefly how to obtain a version of a FCLT when the random walk is replaced by the movement of a particle in a dispersing periodic billiard.
We refer to \cite{Pen09} and \cite{Pen14} for more details on this model.

Let be given a ``billiard table''  in the plane, union of $\Z^2$-periodically distributed obstacles with pairwise disjoint closures. 
We consider a point particle moving in the complementary $Q$ of the billiard table in $\R^2$ with unit speed and elastic reflection off the obstacles. 
By sampling the flow at the successive times of impact with the obstacles, we obtain a Poincar\'e's section of the billiard flow, the billiard transformation.

We assume that the obstacles are strictly convex with pairwise disjoint closures and boundaries of class $C^{r+1}$ with curvature $>0$
(Sinai's billiard or Lorentz's process). Moreover we make the hypothesis of finite horizon (the time between two subsequent reflections is uniformly bounded).  

Suppose that to each obstacle is associated a real random variable with zero expectation, positive and finite variance, 
independent of the motion of the particle and that the family of these r.v.s is i.i.d.

Like in an infinite ``pinball'' with random gain, at each collision with an obstacle, the particle wins the amount given by the random variable associated with the obstacle 
which is met. Let $W_n$ be the total amount won by the particle after $n$ reflections. An annealed FCLT for $W_n$ has been shown by F. P\`ene (\cite{Pen09}): 
there exists $\beta_0 > 0$ such that $\displaystyle {W_{[nt]} \over \beta_0 \, n \lg n}$ converges weakly to the standard Wiener process.

In order to prove a quenched version, we use \cite[Proposition 7]{Pen09}, in place of (\ref{nVarep1}) for the r.w., and \cite[Corollary 4]{Pen14} 
(the main and most difficult step), which gives for the self-intersections of the billiard transformation a law of large numbers replacing (\ref{equivAn20}). 
Then, by 1) in Remarks \ref{general_0} and by the preceding method for the r.w. in random sceneries, we obtain the quenched version of the FCLT for this model.

\section{\bf Cumulants and CLT} \label{Sectcumulant}

For the models of random fields in Sections \ref{movAverag}, \ref{algMod1}, \ref{expMixSect}, we need to recall some tools, 
in particular the method of cumulants which can be used to prove a CLT for dynamical systems satisfying a mixing property of all orders. 

In 1960, Leonov (\cite{Leo60a}, \cite{Leo60b}) applied it to a single algebraic endomorphism of a compact abelian group. 
In \cite{CohCo13}, \cite{CohCo17}, \cite{CohCo16}, it was applied to multidimensional actions given by algebraic endomorphisms
in the connected case and in some non connected cases. 
Recently the method of cumulants has been used in \cite{BjGo20} to prove a CLT for multiple mixing actions with exponential rate.
Using the tightness criterium given in Section \ref{SectMoricz}, we will obtain a functional version in these examples.

\subsection{\bf Moments and cumulants} 

\

For $r \geq 1$, let $X_1, ...,X_r$ be $r$ real centered bounded random variables. We denote by $J_r$ the set $ \{1, ..., r\}$. 
For any subset $I = \{i_1, ..., i_p\} \subset J_r$, we put $m(I) := \E(X_{i_1} ... X_{i_p})$. The {\it cumulant} of order $r$ is
\begin{eqnarray}
C(X_1, ... , X_r) &=&  \sum_{Q = \{I_1, ..., I_p\} \in \Cal Q} (-1)^{p-1} (p - 1)! \ m(I_1) \cdots m(I_p). \label{formCumul0}
\end{eqnarray}
Putting $s(I) := C(X_{i_1}, ... , X_{i_p})$ for $I = \{i_1, ..., i_p\}$, we have
\begin{eqnarray}
\E(X_{1} \cdots X_{r}) &=& m(J_r) = \sum_{Q = \{I_1, ..., I_p\} \in \Cal Q} s(I_1) \cdots s(I_p), \label{cumFormu1}
\end{eqnarray}
where in both formulas, $\Cal Q$ is the set of partitions $Q = \{I_1, I_2, ..., I_p\}$ of $J_r$ into $p \leq r$ nonempty subsets, with $p$ varying from 1 to $r$.

For a single random variable $Y$, the cumulant of order $r$ is defined by $C^{(r)}(Y) := C((Y, ..., Y)_r)$, where $(Y, ..., Y)_r$ is the vector with $r$ components
equal to $Y$.  If $Y$ is centered, we have $\|Y\|_2^2 = C^{(2)}(Y)$ and
\begin{eqnarray}
\E(Y^4) = 3 \E(Y^2)^2 + C^{(4)}(Y). \label{moment4}
\end{eqnarray}
In the next sections, we are going to consider random fields obtained by a measure preserving action $(T^h, h \in \cal H)$ of a group $\cal H$
on a probability space $(E, \cal B, \mu)$. Let $f$ be a measurable bounded centered function on $(E, \mu)$.
For $\underline H =(h_1,..., h_r)$ in $\cal H^r$, we can apply the definition of moments and cumulants to $(T^{h_1}f,..., T^{h_r}f)$.}
\begin{nota} {\rm For the purpose of Section \ref{expMixSect}, we introduce some notations as in \cite{BjGo20}. Let $\underline H$ be in $\cal H^r$ and $I, J $ non empty subsets of $J_r$. We set 
$$d^r(\underline H) := \max_{i,j} d(h_i, h_j), \ d^I(\underline H) := \max_{i, j \in I} d(h_i, h_j), \ d_{I, J}(\underline H) := \min_{i \in I, \, j \in J} d(h_i, h_j)$$
Let $Q$ be a partition of $J_r$, with $|Q| \geq 2$. We set
\begin{eqnarray*}
&&d^Q(\underline H):= \max_{I \in Q} d^I(\underline H), \ d_Q(\underline H):= \min_{I \not = J, \, I,J \in Q} d_{I, J}(\underline H),
\end{eqnarray*}
and, for $0 \leq \alpha < \beta$, 
\begin{eqnarray*}
&&\Delta(\beta) :=\{\underline H \in \cal H^r : d^r(\underline H) \leq \beta\}, 
\ \Delta_Q(\alpha, \beta) :=\{\underline H \in \cal H^r : d^Q(\underline H) \leq \alpha \text{ and } d_Q(\underline H) > \beta\}.
\end{eqnarray*}
}\end{nota}
The elements of a configuration $\underline H = (h_1, h_2, ..., h_r)$ in $\Delta(\beta)$ can be viewed as `clustered', since
$ (h_1, h_2, ..., h_r) \in \Delta(\beta)$ implies $(h_2, ..., h_r) \subset B_d(h_1, \beta)^{r-1}$
 The configurations in $\Delta_Q(\alpha, \beta)$ for some partition $Q$ and $0 < \alpha < \beta$ are made of `well-separated' clusters. 

\vskip 3mm
{\bf Moment of order 4}
 
The moment of order 4 plays a special role in the proof of tightness.

Let $J = [b, b+k]$ be an interval. 
We will bound the moment of order 4 of $Y =\sum_{j \in J} T^{Z_j(\omega)} f$, by using (\ref{moment4}) and by bounding the cumulant of order 4: 
$C^{(4)}(\sum_{j \in J} T^{Z_j(\omega)} f)$. We have:
$$|C^{(4)}(\sum_{j \in J} T^{Z_j(\omega)} f)| 
\leq H(\omega, b , k) := \sum_{s, t, v, w \in J} |C(T^{Z_s(\omega)} f, T^{Z_t(\omega)} f, T^{Z_v(\omega)} f, T^{Z_w(\omega)} f)|.$$
Observe that $H$ is super-additive in the sense of the definition given later in Section \ref{SectMoricz}.

\vskip 3mm
\goodbreak 
{\it Partitions of $\{1, 2, 3, 4\}$}

In Formula (\ref{formCumul0}) of cumulants, the contribution of a partition such that one of its atoms is a singleton is 0, and so does not appear. 
The partitions of $\{1, 2, 3, 4\}$ without atoms reduced to a singleton are $\{\{1, 2, 3, 4\}\}$, $\{\{1, 2\},  \{3, 4\}\}$,  $\{\{1, 3\},  \{2, 4\}\}$, 
 $\{\{1, 4\},  \{2, 3\}\}$. 

The cumulant of order 4 of $(T^{h_i} f, T^{h_k} f, T^{h_\ell} f, T^{h_r} f)$ reads:
\begin{eqnarray}
&&C(T^{h_i} f, \,T^{h_k} f, \,T^{h_\ell} f, \, T^{h_r} f) =  \E(T^{h_i} f \, T^{h_k} f \, T^{h_\ell} f \, T^{h_r} f) \label{cum4} \\
 - &&[\E(T^{h_i}  f \, T^{h_k} f) \, \E(T^{h_\ell} f \, T^{h_r} f) + \E(T^{h_i} f \, T^{h_\ell} f) \, \E(T^{h_k} f \, T^{h_r} f) 
 +\E(T^{h_i}  f \, T^{h_r} f) \, \E(T^{h_k}  f \, T^{h_\ell}  f)]. \nonumber
\end{eqnarray}
 
{\it Well separated configurations}

The following proposition is a key step in the proof of the CLT shown in \cite{BjGo20}.
\begin{prop} \label{SepClust} (cf. Proposition 6.2 in \cite{BjGo20}) For every sequence $(\beta_j)_{j=0, ..., r}$ such that
\begin{eqnarray*}
\beta_0  = 0 < \beta_1 < 3 \beta_1 < \beta_2  < ... < \beta_{r-1} < 3 \beta_{r-1}  <  \beta_r, 
\end{eqnarray*}
we have
\begin{eqnarray}
\cal H^r = \Delta(\beta_r) \cup \bigl(\bigcup_{j=0}^{r-1} \bigcup_{|Q|\geq 2} \, \Delta_Q(3 \beta_j, \beta_{j+1}) \bigr). \label{unionDelta}
\end{eqnarray}
\end{prop}
Below, as an illustration, we give a proof for $r =4$. 
\vskip 3mm
\proof ($r= 4$) Let us consider in general a metric space $\Cal E$ with a distance $d$. By applying the result to $\Cal E =\cal H$, (\ref{unionDelta}) will follows.

We represent an element of  $\Cal E$ by $\bullet$ and, given two elements $a, b \in \Cal E$, 
draw $\bullet \bullet$ or $\bullet \ \bullet$ depending on whether they are close or far from each other.

Given a set with 4 elements in $\Cal E$, $\underline H = \{t, u, v, w\}$, we will show that (up to a permutation) $\underline H$ belongs to one of the configuration types 
$S_0 =  \bullet \bullet \bullet \bullet$, $S_1 = \bullet \bullet \bullet  \ \ \ \bullet$, $S_2 =  \bullet \bullet \ \ \ \bullet \bullet$,
$S_3 = \bullet \bullet \ \ \  \bullet \ \ \ \bullet$, or $S_4 =  \bullet \ \ \ \bullet \ \ \ \bullet \ \ \ \bullet$.

The configurations of type $S_0$ are the `clustered' configurations. This is quantified by saying that these configurations are in $\Delta(\beta)$ for some $\beta > 0$.
The configurations of type $S_4$ are the configurations with pairwise distant elements. They correspond to $\Delta_Q(0, \beta)$ for some $\beta$.

We show that $\underline H$ is either of 
\hfill \break of type $S_0 = \bullet \bullet \bullet \bullet$, with $d(x,y) < \beta_4$, for every $x, y$ in $\underline H$,
 \hfill \break or of type $S_4 = \bullet \ \ \ \bullet \ \ \ \bullet \ \ \ \bullet$, with $d(x,y) \geq \beta_1$, for every $x \not = y$ in $\underline H$,
\hfill \break or of one of the types $S_1$, $S_2$, $S_3$. 

If $\underline H$ is not of type $S_0$ or $S_4$, there are distinct elements, let call them $t, u, w$, such that $d(t, u) < \beta_1$ and $d(t, w) > \beta_4$, which implies
\begin{eqnarray}
d(t, u) < \beta_1, \, d(t, w) > \beta_4 \geq \beta_3, \, d(u, w) > \beta_4 - \beta_1 \geq \beta_3. \label{ineq0}
\end{eqnarray}

\goodbreak
One of the following cases occurs: 

1) (type $S_1 =  \bullet \bullet \bullet  \ \ \ \bullet$) $d(v, t) < 2 \beta_2$, which implies:
\hfill \break $d(v, u) < 2 \beta_2 + \beta_1 < 3 \beta_2$, $d(v, w) > \beta_4 - 2 \beta_2 > 3 \beta_3 - 2 \beta_2 > \beta_3$.

In this case, for the partition $Q = \{\{t, u, v\}, \{w\}\}$, we get $d^Q(\underline H)  < 3 \beta_2$ and $d_Q(\underline H) > \beta_3$.

2) $d(v, t) > 2 \beta_2$, which implies: $d(v, u) > 2 \beta_2 - \beta_1 > \beta_2$,

2a) (type $S_2 =  \bullet \bullet \ \ \ \bullet \bullet$) $d(v, w) < \beta_2$, which implies:
\hfill \break $d(v, t) > d(t, w) - d(v, w) > \beta_4 - \beta_2 > 3 \beta_3 - \beta_2 > \beta_3$,
\hfill \break $d(v, t) > d(t, w) - d(v, w) > \beta_4 - \beta_1 - \beta_2 > 3 \beta_3 - \beta_1 - \beta_2 > \beta_3$.
\hfill \break For $Q = \{\{t, u\}, \{v, w\}\}$, we get $d^Q(\underline H)  < 3 \beta_2$ and $d_Q(\underline H) > \beta_3$.

2b) (type $S_3 = \bullet \bullet \ \ \  \bullet \ \ \ \bullet$) $d(v,w) > \beta_2$ : 
\hfill \break For $Q = \{\{t, u\},\{v\}, \{w\}\}$, we get $d^Q(\underline H)  < 3 \beta_1$ and $d_Q(\underline H) > \beta_2$. \eop

\vskip 3mm
\subsection{\bf A sufficient condition for the CLT}

 \ 

Let us recall a criterium in terms of cumulants for the CLT (cf. \cite[Th. 7]{Leo60b}, \cite[Th. 6.2]{CohCo17}).
It is convenient to formulate the criterium by using a summation sequence, $w= (w_n)_{n \geq 1}$,
i.e., for each $n$ a function $\el \in \Z^d \to w_n(\el) \in \R$, with $0 < \sum_{\el \in \Z^d} |w_n(\el)| < +\infty$. 

The associated normalized non-negative kernel is
$\displaystyle \tilde K(w_n)(\t) = {|\sum_{\el \in \Z^d} w_n(\el) \, e^{2\pi i \langle \el, \t \rangle}|^2 \over \sum_{\el \in \Z^d} |w_n(\el)|^2}, \ \t \in \T^d$.

As for the summation along a random walk which is a special case, we say that the summation is $\xi$-{\it regular} for a probability measure $\xi$ on $\T^d$, 
if the normalised kernel $(\tilde K(w_n)_{n \geq 1})$ converges weakly to $\xi$, i.e., 
\begin{eqnarray}
\lim_{n \to \infty}\int_{\T^d} \tilde K(w_n) \,\varphi \, d\t  = \xi(\varphi), \text{ for every continuous function } \varphi \text{ on } \T^d. \label{regWn}
\end{eqnarray}
This implies that for $f$, under Condition (\ref{sumcorr1}), the asymptotic variance for the sums is
\begin{eqnarray*}
\sigma_w^2(f) := \lim_n {\|\sum_{\el \in \Z^d} w_n(\el) \, T^\el f\|_2^2 \over \sum_{\el \in \Z^d} |w_n(\el)|^2}  = \xi(\varphi_f). \label{asymptVar0}
\end{eqnarray*}
\begin{thm} \label{Leonv0} If $(w_n)_{n \geq 1}$ is a summation sequence on $\Z^d$ such that (\ref{regWn}) holds for a measure $\xi$ on $\T^d$, the condition
\begin{eqnarray}
\sum_{(\el_1, ..., \el_r) \, \in (\Z^d)^r} \, w_n(\el_1) ... w_n(\el_r) \, C(T^{\el_1} f, ... , T^{\el_r} f) = o \bigl((\sum_{\el \in \Z^d} \, w_n^2(\el))^{\frac{r}2}\bigr),
\ \forall r \geq 3, \label{smallCumul0}\\
\text{ implies } \ \  \bigl(\sum_{\el \in \Z^d} \, w_n^2(\el)\bigr)^{-\frac12} \, \sum_{\el \in \Z^d} w_n(\el) f(T^\el .) {\underset{n \to \infty}
\Longrightarrow } \Cal N(0,\xi(\varphi_f)). \label{cvgce1}
\end{eqnarray}
\end{thm}

The following result (cf. \cite[Lemma 6.6]{CohCo17}) shows that mixing of all orders implies the asymptotic nullity of the cumulants.
\begin{proposition} \label{skTozeroLem} Let $(T^\el, \el \in \Z^d)$ be a $\Z^d$-measure preserving action on a probability space $(E, \mu)$. 
If it is mixing of order $r \geq 2$, then, for any $f \in L_0^\infty(X)$, 
\begin{eqnarray}
\underset{\max_{i \not = j} \|\el_i - \el_j\| \to \infty} \lim C(T^{\el_1} f, ... , T^{\el_r} f) = 0. \label{lem66}
\end{eqnarray}
\end{proposition}
Remark that (\ref{lem66}) does not give the quantitative estimate needed in (\ref{smallCumul0}). 
Nevertheless, in Section \ref{algMod1}, (\ref{lem66}) will be sufficient  for an action 
by automorphisms of a connected compact abelian group (in particular of a torus) which is mixing, when $f$ is a trigonometric polynomial. 
For general exponentially mixing actions, a quantitative formulation is needed as in \cite{BjGo20}, using Proposition \ref{SepClust}.

\vskip 2mm
{\it Array of sequences and finite dimensional distributions}

For $s \geq 1$ and $j =1, ..., s$, let $(w_{n, j}, n \geq 1)$ be $s$ summation sequences, satisfying (\ref{regWn}) with respectively $\xi = \xi_j$, 
where the $\xi_j$'s are probability measures on $\T^ d$.

Using Theorem \ref{Leonv0}, we are going to deduce from the following two conditions the asymptotic normality (after normalization) of the vectorial process 
$\bigl(\sum_{\el \in \Z^d} \, w_{n,1}(\el) \, T^\el f, ...,  \sum_{\el \in \Z^d} \, w_{n,s}(\el) \, T^\el f \bigr)$:

- asymptotic orthogonality:
\begin{eqnarray}
\int_{\T^d} (\sum_{\el \in \Z^d} w_{n,j}(\el) \, e^{2\pi i \langle \el, \t\rangle}) \, (\sum_{\el \in \Z^d} w_{n,j'}(\el) \, e^{-2\pi i \langle \el, \t\rangle}) 
\, e^{-2\pi i \langle \p, \t\rangle} \, d\t \nonumber \\
\ \ = o\bigl(\sum_{\el \in \Z^d} (w_{n, j}(\el))^2 + \sum_{\el \in \Z^d} (w_{n, j'}(\el))^2\bigr), \forall j \not = j', \forall \p \in \Z^d. \label{ortho1}
\end{eqnarray}
\hfill \break - convergence to 0 of the normalized cumulants of order $\geq 3$:
\begin{eqnarray}
&&\sum_{(\el_1, ..., \el_r) \, \in (\Z^d)^r} \, w_{n, i_1}(\el_1) ... w_{n, i_r}(\el_r) \, C(T^{\el_1} f, ... , T^{\el_r} f) \nonumber \\
&&\ \ \ = o(\sum_{\el \in \Z^d} \, [\sum_{j=1}^s \, (w_{n, j}(\el))^2])^{r/2}, \, \forall (i_1, ..., i_r) \in \{1, ..., s\}^r, \, \forall r \geq 3. \label{negligCum3}
\end{eqnarray}

\begin{proposition} \label{CltSum1} Under Conditions (\ref{ortho1}) and  (\ref{negligCum3}), the vectorial process 
$$\bigl({\sum_{\el \in \Z^d} \, w_{n, 1}(\el) \, T^\el f \over (\sum_{\el \in \Z^d} w_{n, 1}(\el)^2)^{\frac12}}, ..., 
{\sum_{\el \in \Z^d} \, w_{n, s}(\el) \, T^\el f \over (\sum_{\el \in \Z^d} w_{n, s}(\el)^2)^{\frac12}}\bigr)_{n \geq 1}$$
is asymptotically distributed as ${\cal N}(0, J_s)$, where $J_s$ is the $s$-dimensional diagonal matrix with diagonal $(\xi_j(\varphi_f), j= 1, ..., s)$.
\end{proposition}
\proof The hypothesis (\ref{ortho1}) implies, for $s$ non zero real parameters $a_1, ..., a_s$:
\begin{eqnarray}
(\sum_j a_j^2 \sum_{\el \in \Z^d} (w_{n, j}(\el))^2)^{-1} |\sum_{\el \in \Z^d} \sum_j a_j w_{n, j}(\el) \, e^{2\pi i \langle \el, \t\rangle}|^2 
\, \overset {\text{weakly}} {\underset {n \to\infty} \longrightarrow} \, (\sum_j a_j^2)^{-1} \, \sum_j a_j^2 \, \xi_j. \label{sumVar1}
\end{eqnarray}
Putting $w_n^{a_1, ..., a_s}(\el) = a_1 w_{n, 1}(\el) + ... + a_s w_{n, s}(\el)$, 
by the Cram\'er-Wold theorem, to conclude it suffices to show 
\begin{eqnarray}
{\sum_{\el \in \Z^d} \, w_n^{a_1, ..., a_s}(\el) T^\el f \over (a_1^2 \sum_{\el \in \Z^d} (w_{n, 1}(\el))^2 + ... + a_s^2  \sum_{\el \in \Z^d} \, (w_{n, s}(\el))^2)^{\frac12}}
\Longrightarrow \Cal N(0, \sum_{j=1}^s a_j^2 \,\xi_j(\varphi_f) / \sum_{j=1}^s a_i^2). \label{TCLFd1}
\end{eqnarray}
By (\ref{negligCum3}), the sum $\displaystyle \sum_{i_1, ..., i_r \in \{1, ..., s\}^r} \, \sum_{(\el_1, ..., \el_r) \, \in (\Z^d)^r} \, w_{n, i_1}(\el_1) ... w_{n, i_r}(\el_r) 
\, C(T^{\el_1} f, ... , T^{\el_r} f)$ satisfies (\ref{smallCumul0}) and the result follows from Theorem \ref{Leonv0}. \eop

\vskip 3mm
We will use also the following lemma. Let $(w_n)_{n \geq 1}$ be a summation sequence on $\Z^d$ such that (\ref{regWn}) holds 
for a measure $\xi$ on $\T^d$. For $f \in L^2(\mu)$, we put $\sigma_n(f) := \|\sum_\el w_n(\el) \, T^\el f \|_2$.
We can suppose $\xi(\varphi_f) > 0$, since otherwise the limiting distribution is $\delta_0$. 
\begin{lem} \label{tclRegKer} Let $f, f_k, k \geq 1$ be in $ L^2(\mu)$ and satisfying (\ref{sumcorr1}) such that $\|\varphi_{f - f_k}\|_\infty \to 0$. Then
\begin{eqnarray*}
&&{\sum_{\el \in \Z^d} w_n(\el) \, T^\el f_k \over \sigma_n(f_k)} \  {\underset{n \to \infty} \Longrightarrow} \ \Cal N(0,1), \, \forall k  \geq 1, 
\text{ implies } \, {\sum_{\el \in \Z^d} w_n(\el) \, T^\el f \over \sigma_n(f)} \ {\underset{n \to \infty} \Longrightarrow} \ \Cal N(0,1). 
\end{eqnarray*}
\end{lem}
\proof Let $(\varepsilon_k)$ be a sequence of positive numbers tending to 0, such that $\|\varphi_{f - f_k}\|_\infty \leq \varepsilon_k$. 
Let us consider the processes defined respectively by
\begin{eqnarray*}
U_n^{k} := (\sum_{\el \in \Z^d} \, w_n^2(\el))^{-\frac12} \, \sum_{\el \in \Z^d} w_n(\el) \, T^\el f_k, 
\ U_n := (\sum_{\el \in \Z^d} \, w_n^2(\el))^{-\frac12} \, \sum_{\el \in \Z^d} w_n(\el) \, T^\el f.
\end{eqnarray*}
By (\ref{asymptVar0}) we have:
$$(\sum_{\el \in \Z^d} \, w_n^2(\el))^{-1} \, \|\sum_{\el \in \Z^d} w_n(\el) \, T^\el f \|_2^2 
= \int_{\T^d} \, \tilde w_n \, \varphi_f \, dt \underset{n \to \infty} \to \xi(\varphi_f).$$
Since $\xi(\varphi_{f - f_k}) \to 0$, we have $\xi(\varphi_{f_k}) \not = 0$ for $k$ big enough.
Indeed, if $\xi$ is a probability measure on $\T^d$, $f \to (\xi(\varphi_f))^\frac12$ satisfies the triangular inequality.

The hypotheses imply ``${U_n^{k} \ {\underset{n \to \infty} \Longrightarrow} \Cal N(0,\xi(\varphi_{f_k}))}$'' for every $k$. Moreover, since
\begin{eqnarray*}
\lim_n \int |U_n^k - U_n|_2^2 \ d\mu &=& \lim_n \int_{\T^d} \, \tilde w_n \,\varphi_{f-f_k} \, d \t = \xi(\varphi_{f - f_k}) \leq \varepsilon_k,
\end{eqnarray*}
we have $\limsup_n \mu[|U_n^k - U_n| > \delta] \leq \delta^{-2} \limsup_n \int |U_n^k - U_n|_2^2 \ d\mu \underset{k \to \infty} \to 0$, for every $\delta > 0$.

Therefore, using \cite[Theorem 3.2]{Bi99}, the conclusion ``$U_n \ {\underset{n \to \infty} \Longrightarrow } \Cal N(0, \xi(\varphi_{f}))$'' follows. \eop

\section{\bf Moving averages of i.i.d. random variables} \label{movAverag}

Let $(X_\el)_{\el \in \Z^2}$ be a r.f. of centered i.i.d. real random variables such that $\|X_\0\|_2 = 1$. Let $(a_\q)_{\q\in \Z^2}$ be an array of real numbers
such that $\sum_{\q \in \Z^2} |a_\q| < \infty$ and let  $\Xi = (\Xi_\el)_{\el \in \Z^2}$ be the random field defined by $\Xi_\el(x) = \sum_{\q \in \Z^2} a_\q X_{\el - \q}(x)$. 

The correlation is 
$\displaystyle \widehat \varphi_\Xi(\el) = 
\langle \sum_{\q \in \Z^2} a_\q X_{\el - \q}, \sum_{\q' \in \Z^2} a_{\q'} X_{- \q'} \rangle = \sum_{\q \in \Z^2} a_\q \, a_{\q - \el}$.  
We have 
$$\sum_\el |\widehat \varphi_\Xi(\el)| \leq \sum_\el \sum_{\q \in \Z^2} |a_\q| \, |a_{\q - \el}| = (\sum_{\q \in \Z^2} |a_\q|)^2< +\infty.$$
The continuous spectral density of the process is $\varphi_\Xi(\t)  = |\sum_\q a_\q e^{2\pi i \langle\q, \t\rangle}|^2$.
The asymptotic variance for the summation along the r.w. (with the normalisation by $C_0 n \ln n$) is $(\sum_{\q \in \Z^2} a_\q)^2$
that we suppose $\not = 0$.

Using the method of associated r.v.s we obtain a quenched FCLT for $S_{\lfloor nt \rfloor}^{\omega, \Xi}$ (cf. Notation (\ref{defsumRW})): 
\begin{thm} \label{FCLTindFil} The process $\displaystyle \bigl({S_{\lfloor nt \rfloor}^{\omega, \Xi}(x) \over\sqrt{ n \log n}}\bigr)_{t \geq 0}$
satisfies a quenched FCLT with asymptotic variance $\sigma^2 = |\sum_{\q \in \Z^2} a_\q |^{2} (\pi \sqrt{\det \Sigma})^{-1}.$
\end{thm}
\proof 1) {\it Convergence of the finite dimensional distributions}

a) First we assume that the random variables are bounded. 
Moreover let us consider first a finite sum $F = \sum_{\s \in S} a_\s X_\s$, where $S$ is a finite subset of $\Z^2$. 
The case of the series, $\Xi_\0 = \sum_{\s \in \Z^2} a_\s X_\s$,  will follow by an approximation argument.

As we have seen, for $Y_n(\omega, x, t) = {S_{[n t]}^{\omega, X}(x) \over \sqrt{C_0 n \log n}}$, we have to show:
$$\forall \, 0=t_0 < t_1 <... < t_r =1, \ (Y_{n}(t_1), ..., Y_{n}(t_r)) \underset {n \to \infty} \Longrightarrow (W_{t_1}, ..., W_{t_r}).$$
For it, we use Proposition \ref{CltSum1}. Condition (\ref{ortho1}) follows from  Lemma \ref{VarrwLem1}.  Let us check (\ref{negligCum3}).

There is $M$ such that the cumulant $C(T^{\el_1} F,..., T^{\el_r} F) = 0$, if $\max_{i,j} \|\el_i - \el_j\| > M$, because if $M$ is big enough, there is a random variable
$T^{\el_{i}} F$ which is independent from $\sigma$-algebra generated by the others in the collection $T^{\el_1} F, ..., T^{\el_r} F$ (by finiteness of $S$). 

Let $\displaystyle w_n^{u}(\omega, \el) := \sum_{j=n t_u}^{nt_{u+1}} 1_{Z_j=\el}$, for $u= 0, ..., r-1$. Then we have $w_n^{u}(\omega, \el) \leq w_n(\omega, \el)$
and, since $\displaystyle \sup_{\el_1, ..., \el_r} |C(T^{\el_1} F, T^{\el_2} F, ..., T^{\el_r} F)| < \infty$, 
\begin{eqnarray*}
&& |\sum_{\max_{i,j} \|\el_i - \el_j\| \leq M} C(T^{\el_1} F, T^{\el_2} F, ..., T^{\el_r} F)| \, 
w_n^{i_1}(\omega, \el_1) \, w_n^{i_2}(\omega, \el_2) ... w_n^{i_r}(\omega, \el_r) \\
&&\ \ \leq \sum_{\el} \sum_{\|\j_2\|, ..., \|\j_r\| \leq M, \, \j_1 =\0} |C(T^{\el} F, T^{\el + \j_2} F, ..., T^{\el + \j_r} F)| \, \prod_{k=1}^r w_n^{i_k}(\omega, \el+\j_k)\\ 
&&\ \leq C \sum_\el \sum_{\|\j _2\|, ..., \|\j_r\| \leq M, \, \j_1 =\0} \ \prod_{k=1}^r w_n^{i_k}(\omega, \el+j_k)
\ \leq C \sum_\el \sum_{\|\j _2\|, ..., \|\j_r\| \leq M, \, \j_1 =\0} \ \prod_{k=1}^r w_n(\omega, \el+j_k).
\end{eqnarray*}
The right hand side is less than a finite sum of sums of the form $\sum_{\el \in \Z^d} \prod_{k= 1}^r w_n(\omega, \el + {\j}_k)$ with $\{\j_1, ..., \j_r\} \in \Z^2$.

By (\ref{nVarep1}), for every $\varepsilon > 0$, there is $C_\varepsilon(\omega)$ a.e. finite such that 
$\sup_\el w_n(\omega, \el) \leq C_\varepsilon(\omega) \, n^\varepsilon$. For $r \geq 3$, take $\varepsilon < {r-2 \over 2(r-1)}$. We have then
$\sum_{\el \in \Z^d} \prod_{k= 1}^r w_n(\omega, \el + {\j}_k) \leq  C_\varepsilon(\omega)^{r-1} \, n^{\varepsilon(r-1)}\, n = o(n^{r/2})$
and (\ref{negligCum3}) is satisfied. 

Using Lemma \ref{tclRegKer}, the result can be extended to a general sum $\sum_{\s \in S} a_\s X_\s$, such that $\sum_{\s \in S} |a_\s| < \infty$.

b) Now if we assume only the condition $\|X_\0\|_2 < \infty$, we use a truncation argument and apply again Lemma \ref{tclRegKer}.

2) {\it Tightness} Let $a_\q^+= \max (a_\q, 0)$, $a_\q^-= \max (-a_\q, 0)$.
Observe that the random variables $\sum_{\q \in \Z^2} a_\q^+ X_{\el - \q}(x) = \sum_{\q \in \Z^2} a_{\q+\el}^+ X_{- \q}(x)$, for $\el \in \Z^2$, are associated, 
as well as $\sum_{\q \in \Z^2} a_\q^- X_{\el - \q}(x)$, for $\el \in \Z^2$.

Therefore tightness can be proved separately for both processes. The proof is like the proof of tightness in Theorem \ref{FCLTind1}. \eop

\section{\bf Tightness and 4th-moment} \label{SectMoricz}

In this section, we show a criterium of tightness based on the 4th-moment.  

Let $n \geq 1$. We say that a nonnegative function $G_0= (G_0(b, k))$, defined for $b, k$ such that $0 \leq b \leq b+k \leq n$, is super-additive if $G_0(b, 0) = 0$ and 
\begin{eqnarray}
G_0(b,k) + G_0(b+k, \ell) \leq G_0(b,k+\ell) , \, \forall b \geq 0, \forall k, \ell \geq 1 \text{ such that } b+k+\ell \leq n. \label{supadd1}
\end{eqnarray} 

Let $(W_k)$ be a sequence of real or complex random variables on a probability space $(E, \mu)$. We set
$$S_{b,k} = \sum_{r=b+1}^{b+k} W_r, \, M_{b,n} = \max_{1 \leq k \leq n} |S_{b,k}|.$$
The following result is adapted from \cite{Mo76}:
\begin{thm} (F. M\'oricz) \label{thmMoricz1} Let $n \geq 1$. Suppose that there exists $G_0$ a  super-additive function such that
\begin{eqnarray}
\E_\mu(|S_{b, k}|^4) \leq G_0^2(b, k), \, \forall \, b, k \text{ such that } 0 \leq b \leq b+k \leq n. \label{hypoG1}
\end{eqnarray} 
Then, with the constant $C_{max} = (1 - 2^{-\frac14})^{-4}$,
\begin{eqnarray}
\E_\mu(|M_{b,n}|^4) \leq C_{max} \, G_0^2(b,n),  \, \forall b \leq n. \label{toprove1}
\end{eqnarray} 
\end{thm}
Let $X = (X_\el)_{\el \in \Z^2}$ be a strictly stationary real r. f. on a probability space $(E, \mu)$, where the $X_\el$'s have zero mean and finite second moment.
Setting $S_J^{\omega}(x) =\sum_{i \in J} X_{Z_i(\omega)}(x)$ if $J$ is an interval, we deduce from (\ref{toprove1}) a criterium for tightness 
adapted to the sums along a random walk.

\begin{proposition} \label{rwFCLTB} Let $G(\omega, ., .)$, $H(\omega, ., .)$ be super-additive functions such that for a parameter $\gamma$ 
and $K_1(\omega), K_2(\omega)$ a.e. finite functions on $(\Omega, \PP)$,
\begin{eqnarray}
G(\omega, b, k) \leq  K_1(\theta^b \omega) \, k \ln k, \ H(\omega, b, k) \leq  K_2(\theta^b \omega) \, k \, (\ln k)^{\gamma}, \ G(\omega, b, k) \geq k. \label{majGH1}
\end{eqnarray}  
Suppose that the r.v.s $X_\el$ are bounded and satisfy
\begin{eqnarray}
&&\E_\mu(|S_{J}^{\omega, X}|^4) \leq G(\omega, b, k)^2 +  n^\frac12 \, (\ln n)^{-(\gamma + 1)}\, H(\omega, b, k), \forall J= [b, b+k] \subset [1, n], 
\text{ for a.e. } \omega. \label{hypoG1b}
\end{eqnarray} 
Then, for every $\varepsilon > 0$, $Y_n(\omega, x, t) = {1 \over \sqrt {n \ln n}} \sum_{j=1}^{[nt]} \, X_{Z_j(\omega)}$ satisfies
\begin{eqnarray}
&& \lim_{ \delta \to 0^+} \limsup_n \mu(x \in E: \, \sup_{|t'-t| \leq \delta}|Y_n(\omega, x, t') - Y_n(\omega, x, t)| \geq \varepsilon) = 0. \label{tight1}
\end{eqnarray}
\end{proposition}
\proof 1) Let $c \geq 0$, $\Delta_n = n^\frac12 (\ln n)^{-2}$, $\nu= \nu_n \geq \Delta_n$,
$L_n = [{\nu_n \over \Delta_n}]$, $\nu' = \nu_n' = [{\nu_n \over \Delta_n}] \Delta_n + \Delta_n + 1$.

The integer $\nu_n$ will be chosen of order $\delta n$. We can write, with the convention that $\sum_{r=0}^{-1} =0$: 
\begin{flalign*}
&\max_{0\leq k \leq \nu}|\sum_{j=0}^{k} X_{Z_{j + c}(\omega)}| \leq \max_{0\leq k \leq \nu'}|\sum_{j=0}^{k} X_{Z_{j + c}(\omega)}|
=\max_{0 \leq u \leq[{\nu\over \Delta_n}], 1 \leq k \leq \Delta_n -1}|\sum_{r=0}^{u-1} 
\sum_{j=r\Delta_n}^{(r+1)\Delta_n -1} X_{Z_{j + c}(\omega)}+ \sum_{j = u \Delta_n}^{u \Delta_n + k - 1} X_{Z_{j + c}(\omega)}|&\\
&\leq\max_{0 < u \leq L_n, \, 1 \leq k \leq \Delta_n -1}|\sum_{r=0}^{u-1} \, \sum_{j=r\Delta_n}^{(r+1)\Delta_n -1} \, X_{Z_{j + c}(\omega)}| 
+ \max_{0 \leq u \leq L_n, \, 1 \leq k \leq \Delta_n -1}\, |\sum_{j = u \Delta_n}^{u \Delta_n + k - 1} \, X_{Z_{j + c}(\omega)}|&\\
&= \max_{0 < u \leq L_n} \, |\sum_{j=0}^{u\Delta_n -1} \, X_{Z_{j + c}(\omega)}| 
+ \max_{0 \leq u \leq L_n, \, 1 \leq k \leq \Delta_n -1}\, |\sum_{j = u \Delta_n}^{u \Delta_n + k - 1} \, X_{Z_{j + c}(\omega)}| = \hat A_n + \tilde A_n.&
\end{flalign*}

With $\hat A_n$ and $\tilde A_n$ respectively the first and the second term above, this implies 
\begin{eqnarray}
&&\mu(\max_{0\leq k \leq \nu}|\sum_{j=0}^{k} \, X_{Z_{j + c}(\omega)}|  \geq \varepsilon \sqrt {n \ln n}) 
\leq \mu(\hat A_n \geq \frac12\varepsilon \sqrt {n \ln n}) + \mu(\tilde A_n \geq \frac12\varepsilon \sqrt {n \ln n}). \label{discr1}
\end{eqnarray}
For $\tilde A_n$, since the $X_\el$'s are bounded (uniformly in $\el$ by stationarity), by the choice of $\Delta_n$
there is $N_1(\varepsilon, \delta)$ such that $\mu(\tilde A_n \geq \frac12\varepsilon \sqrt {n \ln n}) = 0$, for $n \geq N_1(\varepsilon, \delta)$.

For $\hat A_n$ we will apply Theorem \ref{thmMoricz1} to $W_r = \sum_{j=r\Delta_n}^{(r+1)\Delta_n -1} \, X_{Z_{c+j}(\omega)}$, with
\begin{eqnarray}
G_0(b, k) := G(\omega, c+b \Delta_n, k \Delta_n) + (\ln n)^{- \gamma + 1} \, H(\omega, c + b \Delta_n, k \Delta_n), \label{defG0}
\end{eqnarray}
which is super-additive as $G$ and $H$.

Since $G(\omega, c+b \Delta_n, k\Delta_n) \geq G(\omega, c+b \Delta_n, \Delta_n) \geq \Delta_n =  n^\frac12 (\ln n)^{-2}$, we have for $k \geq 1$:
\begin{eqnarray*}
&&G^2(\omega, c+b \Delta_n,  k\Delta_n) + n^\frac12 \, (\ln n)^{-\gamma - 1}\, H(\omega, c+b \Delta_n, k\Delta_n) 
\leq G_0(b, k)^2,
\end{eqnarray*} 
Therefore, by (\ref{hypoG1b}),
\begin{eqnarray*}
&&\E_\mu(|\sum_{r=b+1}^{b+k} \, W_r|^4)
= \E_\mu(|\sum_{j=(b+1)\Delta_n}^{(b+k+1)\Delta_n -1} \, X_{Z_{c+j}(\omega)}|^4) \leq G_0(b, k)^2, \, \forall b \geq 0, \forall k \geq 1,
\end{eqnarray*}
which implies by (\ref{toprove1}) of Theorem \ref{thmMoricz1}: 
\begin{eqnarray*}
&&\E_\mu(\max_{1 \leq k \leq p} |\sum_{j=(b+1)\Delta_n}^{(b+k+1)\Delta_n -1} \, X_{Z_{c+j}(\omega)}|^4)
\leq C_{\max} \, G_0(b, p)^2,  \, \forall b \geq 0, \forall p \geq 1.
\end{eqnarray*} 
Putting $K(\omega) := \max(K_1(\omega), K_2(\omega))$ and using (\ref{majGH1}), we get the bound
\begin{eqnarray}
&&\|\max_{u = 1}^{L_n} |\sum_{j=0}^{u\Delta_n} \, X_{Z_{c+j}(\omega)}| \|_4^4 
\leq C_{max} \, [G(\omega, c, L_n \Delta_n) + (\ln n)^{-\gamma + 1} \, H(\omega, c, L_n \Delta_n)]^2 \nonumber \\
&&\leq  C_{max} \, K(\theta^c \omega)^2 \, [L_n \Delta_n \ln (L_n \Delta_n) + (\ln n)^{-\gamma + 1} \, (L_n \Delta_n) (\ln (L_n \Delta_n))^{\gamma}]^2. \label{upBound1}
\end{eqnarray}

2) For $M > 0$ big enough, the set $\Omega_M := \{\omega: \, K(\omega) \leq M\}$ has a probability $\PP(\Omega_M) \geq \frac12$. 
We apply Lemma \ref{mean1} to $\Omega_M$. Given $\delta > 0$, there is $N_2(\delta)$ such that for $n \geq N_2(\delta)$, 
we can find a sequence $0 = \rho_{1, n} < \rho_{2, n} < ... < \rho_{v, n} \leq n < \rho_{v+1, n}$ of visit times of $\theta^k \omega$ in $\Omega_M$ 
under the iteration of the shift $\theta$, such that $\frac12 \delta n \leq \rho_{i+1, n} - \rho_{i, n} \leq \frac32 \delta n$ and $v < 2/\delta$.
By construction, $K(\theta^{\rho_{i, n}}\omega) \leq M, \forall i$.

With $c= \rho_{i, n}$, $\nu_n=\nu_{i, n} = \rho_{i+1, n} - \rho_{i, n} \leq \frac32 \delta n$, $L_{i, n} = [{\nu_{i, n} \over \Delta_n}]$
(so that $L_{i, n} \Delta_n \leq \frac32 \delta n$), we deduce from the upper bound (\ref{upBound1}) 
(for $n$ big enough and using $0 \leq \ln (\delta n) \leq \ln n$, if $n \geq \delta^{-1}$):
\begin{eqnarray*}
&&\|\max_{u = 1}^{L_{i, n}} |\sum_{j=0}^{u\Delta_n} \, X_{Z_{\rho_{i, n}+j}(\omega)}| \|_4^4 \leq C_{max} M \, [\nu_{i, n} \, \ln \nu_{i, n} 
+ (\ln n)^{-\gamma + 1} \, \nu_{i, n} \, (\ln \nu_{i, n})^{\gamma}]^2\\
&&\leq C_{max} M \, [\frac32\delta n \ln (\delta n) + (\ln n)^{-\gamma + 1} \, \frac32 \delta n (\ln (\delta n))^{\gamma}]^2 \leq C_{max} M \, [3\delta n \ln n]^2.
\end{eqnarray*}
This implies for $\hat A_n$ (cf. (\ref{discr1})), for $i= 1,  ..., v$, for a constant $C$:
\begin{flalign*}
&\mu(\max_{0 \leq u \leq L_{i, n}}|\sum_{j=0}^{u\Delta_n -1} \, X_{Z_{j + \rho_i}(\omega)}| \geq \frac12 \varepsilon \sqrt {n \ln n})
\leq {C_{max} M \, (3 \delta n \ln n)^2 \over (\frac12 \varepsilon \sqrt {n \ln n})^4} \leq C \, \varepsilon^{-4} \, \delta^2.&
\end{flalign*}
Putting $t_i = \rho_{i, n} /n$, we obtain, for $n \geq N(\varepsilon, \delta)$ with $N(\varepsilon, \delta) \geq N_1(\varepsilon, \delta)$ and big enough,
\begin{flalign*}
&\mu (\sup_{|t'-t| \leq \delta}|Y_n(t') - Y_n(t)| \geq 3\varepsilon) \leq \sum_{i=1}^v \mu(\sup_{t_{i-1}\leq s \le t_i}| Y_n(s) - Y_n(t_{i-1})| \geq \varepsilon)
\leq 2 C \, \varepsilon^{-4} \, \delta^2 \, v \leq 2 C \, {\delta \over \varepsilon^4}. \eop&
\end{flalign*}

\begin{rem} \label{remSum1} Let be given for each $s$ in a set of indices $S$ a process $X^s = (X_{s, \el})_{\el \in \Z^2}$ satisfying the hypotheses of the proposition, 
with the same uniform bound and the same $G, H, \gamma$. Then, if $X_\el = \sum_s a_s X_{s, \el}^s$ with $\sum_s |a_s| \leq 1$, 
the r.f. $X = (X_\el)$ satisfies the conditions of the proposition and therefore the conclusion (\ref{tight1}).
This follows from Minkowski inequality:
\begin{eqnarray*}
\|S_{J}^{\omega, X}\|_4^4 &\leq& (\sum_s |a_s| \|S_{J}^{\omega, {X^s}}\|_4)^4 \leq (\sum_s |a_s| [G(\omega, b, k)^2 
+ n^\frac12 \, (\ln n)^{-(\gamma + 1)} \ H(\omega, b, k)]^\frac14)^4\\
&=& (\sum_s |a_s|)^4 [G(\omega, b, k)^2 + n^\frac12 \, (\ln n)^{-(\gamma + 1)} \, H(\omega, b, k)].
\end{eqnarray*}
\end{rem}

\section{\bf Random walks and FCLT for automorphisms of a torus \label{algMod1}} 

We consider now a random field generated by the action of commuting automorphisms on a torus. Let us first present the model.
We will give the details of the proof for $d=2$.

{\it Actions by endomorphisms on a compact abelian group:}
Let $G$ be a compact abelian group with Haar measure $\mu$. The group of characters of $G$ is denoted by $\hat G$ and the set of non trivial characters by $\hat G^*$. 
The Fourier coefficients of a function $f$ in $L^1(G, \mu)$ (denoted also $\hat f(\k)$ when $G$ is a torus) are $c_f(\chi) := \int_G \, \overline \chi \, f \, d\mu$, $\chi \in \hat G$.

Every surjective endomorphism of $G$ defines a measure preserving transformation on $(G, \mu)$ and a dual injective endomorphism on $\hat G$. 

Let $(T_1, ..., T_d)$ be a finite family of $d$ commuting surjective endomorphisms of $G$ and $T^\el = T_1^{\ell_1}... T_1^{\ell_1}$, 
for $\el = (\ell_1, ..., \ell_d) \in \Z^d$. We obtain a $\Z^d$-action $(T^\el, \el \in \Z^d)$ on $G$, which is totally ergodic if and only if the dual action is free. 

Let $AC_0(G)$ denote the space of real functions on $G$ with {\it absolutely convergent Fourier series} and $\mu(f) =0$, endowed with
the norm: $\|f \|_c := \sum_{\chi \in \hat G} |c_f(\chi)| < +\infty$. 

Recall that the action on $G$ is mixing of all orders if it is totally ergodic and $G$ is connected.
\begin{prop} \label{resum} If $f$ is in $AC_0(G)$, the spectral density $\varphi_f$ is continuous on $\T^\rho$ and $\|\varphi_{f}\|_\infty \leq \|f\|_c^2$. 
For every $\varepsilon > 0$ there is a trigonometric polynomial $P$ such that $\|\varphi_{f - P}\|_\infty \leq \varepsilon$.
\end{prop}
\proof \ Since by total ergodicity the characters $T^\el \chi$ for $\el \in \Z^d$ are pairwise distinct, we have 
$$\sum_{\el \in \Z^d} |\langle T^{\el}f, f\rangle| \leq \sum_{\el \in \Z^d} \sum_{\chi \in \hat G^*} |c_f(T^\el \chi)| \, |c_f(\chi)|
\leq \sum_{\chi \in \hat G^*} (\sum_{\el \in \Z^d} |c_f(T^\el \chi)|) \, |c_f(\chi)| \leq (\sum_{\chi \in \hat G^*} |c_f(\chi)|)^2.$$
Therefore, if $f$ is in $AC_0(G)$, then $\sum_{\el \in \Z^d} |\langle T^{\el}f, f\rangle| < \infty$, the spectral density is continuous and 
$\|\varphi_f\|_\infty \leq \varepsilon$. By this inequality, we can take for $P$ the restriction of the Fourier series
of $f$ to a finite set  $\Cal E$ in $\hat G$, where $\cal E$ is such that
$\|\varphi_{f - P}\|_\infty \leq (\sum_{\chi \in \hat G \setminus {{\Cal E}}} |c_f(\chi)|)^2 \leq \varepsilon$. \eop

For compact abelian groups which are connected (cf. \cite{CohCo17}) or which belong to a special family of non connected groups (cf. \cite{CohCo16}),
a CLT has been shown for summation either over sets or along a random walk. 
Our aim is to extend this latter result to a functional CLT at least in the case of automorphisms of a torus.

{\it Matrices and automorphisms of a torus:}
Now we will restrict to the special case of matrices and automorphisms of $G = \T^\rho$, $\rho > 1$.

Every $A$ in the semigroup ${\Cal M}^*(\rho, \Z)$ of non singular $\rho \times \rho$ matrices with coefficients in $\Z$ defines a
surjective endomorphism of $\T^\rho$ and a measure preserving transformation on $(\T^\rho, \mu)$. It defines also a dual
endomorphism of the group of characters $\widehat {\T^\rho}$ identified with $\Z^\rho$ (this is the action by the transposed matrix, but since
we compose commuting matrices, for simplicity we do not write the transposition). The linear operator on $\C^\rho$ defined by $A$ is denoted by $\tilde A$.

When $A$ is in the group $GL(\rho, \Z)$ of matrices
with coefficients in $\Z$ and determinant $\pm 1$, it defines an automorphism of $\T^\rho$.
Recall that the action of $A \in {\Cal M}^*(\rho, \Z)$ on $(\T^\rho, \mu)$ is ergodic if and only if $A$ has no eigenvalue root of unity. 

Here we present the proof for  the case of automorphisms and for $d= 2$ (in the recurrent case for the random walk). 
Let $(A_1, A_2)$ be two commuting matrices in $GL(\rho, \Z)$
and $A^\el = A_1^{\ell_1} \, A_2^{\ell_2}$, for $\el = (\ell_1, \ell_2) \in \Z^2$. It defines a $\Z^2$-action $(A^\el, \el \in \Z^2)$ on $(\T^\rho, \mu)$, 
which is totally ergodic if and only if $A^\el$ has no eigenvalue root of unity for $\el \not = \0$.

Explicit totally ergodic $\Z^2$-actions can be computed (cf. \cite{CohCo17}) like the example below 
(see the book of H. Cohen on computational algebraic number theory \cite{Co93}):

$A_1 = \left(
\begin{matrix} -3 & -3 & 1 \cr 10 & 9 & -3 \cr -30 & -26 & 9 \cr
\end{matrix} \right), \ \ A_2  = \left(
\begin{matrix} 11 & 1 & -1 \cr -10 & -1 & 1 \cr 10 & 2 & -1 \cr
\end{matrix}
\right).$ 

\vskip 3mm
We will need an algebraic result based on the following theorem on S-unit equations (\cite{Schl90}):
\begin{thm} \label{EvScSc} (\cite[Th. 1.1]{EvScSc02}) \ Let $K$ be an algebraically closed field of characteristic 0 and
for $r \geq 2$, let $\Gamma_r$ be a subgroup of the multiplicative group $(K^*)^r$ of finite rank.
For any $(a_1,..., a_r) \in (K^*)^r$, the number of solutions $x = (x_1,... , x_r) \in \Gamma_r$ of the equation
$a_1 x_1 + ... + a_r x_r = 1$ such that no proper subsum of $a_1 x_1 + ... + a_r x_r$ vanishes, is finite.
\end{thm}
\begin{cor} \label{nbsolu1} Suppose that the $\Z^2$-action $(A^\el, \el \in \Z^2)$ is totally ergodic.
The set $F$ of triples $(\el_1, \el_2, \el_3) \in (\Z^2)^3$ for which there is $\gamma \in \Z^2 \setminus \{\0\}$ such that
\begin{eqnarray} 
&&A^{\el_1} \gamma - A^{\el_2} \gamma + A^{\el_3} \gamma - \gamma  = 0, \label{unitEq1}
\end{eqnarray}
without vanishing proper sub-sum, is finite.
\end{cor}
\proof There exists a decomposition of $E = \C^\rho$ into vectorial subspaces $\C^\rho = \oplus_k E_k$ which are simultaneously invariant by $\tilde A_i$, $i=1, 2$, 
and such that  there is a basis $B_k$ in which $\tilde A_i$ restricted to $E_k$ is represented in a triangular form with an eigenvalue 
of $\tilde A_i$ on the diagonal.

This follows from the fact that the commuting matrices $A_i$ have a common non trivial space $W$ of eigenvectors,
and then by an induction on the dimension of the vector space, applying the induction hypothesis to the action of the quotient map of $\tilde A_i$ on $E/W$.

For $\gamma \in \Z^\rho \setminus \{0\}$, there is $k_0$ such that the component $\gamma_0$ of $\gamma$ in $E_{k_0}$ is $\not = 0$.
Let $\delta_0$ be the dimension of $E_{k_0}$. In the basis $B_{k_0} = \{e_{k_0, 1}, ..., e_{k_0, \delta_0}\}$ of $E_{k_0}$, we denote the coordinates 
of $\gamma_0$ by $(\gamma_0^1, ..., \gamma_0^{\delta_0})$. There is $\delta_0' \in \{1, ..., \delta_0\}$ such that
$\gamma_0^i = 0, \, \forall i < \delta_0', \text{ and }  \v_0 := \gamma_0^{\delta_0'} \not = 0$.

Due to the triangular form, for $j = 1, 2$, we have $A_j^\ell \gamma_0 = \alpha_{k_0, j}^\ell \, \gamma_0^{\delta_0'} \, e_{k_0, \delta_0'} + \zeta(j, \ell)$, $\forall \ell \in \Z$, 
where $\alpha_{k_0, j}$ is an eigenvalue of $A_j$ and where $\zeta(j, \ell)$ belongs to the subspace generated by $\{e_{k_0, \delta_0' +1}, ..., e_{k_0, \delta_0}\}$.

Using the notation $\al_u^\el = \alpha_{u, 1}^{\ell^1} \, \alpha_{u, 2}^{\ell^2}$, if $\alpha_{u, 1}$ (resp. $\alpha_{u, 2}$) is an eigenvalue of $A_1$ (resp. $A_2$),
if (\ref{unitEq1}) holds, then $(\al_{k_0}^{\el_1} - \al_{k_0}^{\el_2} + \al_{k_0}^{\el_3}) \v_0  = \v_0$. This equation is still without vanishing proper sub-sum, 
because of the assumption of total ergodicity. By Theorem \ref{EvScSc} applied to the multiplicative group (of finite rank) generated by $\alpha_{k_0, j}$, $j=1, 2$,
the number of solutions of the previous equation is finite. Hence the result, since $k_0$ takes a finite number of values.
\eop

\vskip 3mm
{\bf Random walks and quenched CLT} \label{rwSect}

Our aim is to replace the r.f. of i.i.d. variables $(X_\el, \el \in \Z^2)$ discussed in Section \ref{indepSec} by 
the random field generated by an observable $f$ on a torus $\T^\rho$ under the action of commuting automorphisms. 

More precisely, we consider $\el\mapsto A^{\underline \ell}$ a totally ergodic $\mathbb Z^2$-action by algebraic automorphisms of $\mathbb T^\rho$, $\rho>1$, 
defined by commuting $\rho\times \rho$ matrices $A_1, A_2$ with integer entries, determinant $\pm 1$ 
such that the eigenvalues of $A^{\el} = A_1^{\ell_1} A_2^{\ell_2}$ are $\not = 1$, if $\el = (\ell_1, \ell_2) \not = (0, 0$).

The composition with a function $f$ defined on $\T^\rho$ is denoted by $A^\el f$ as well as $T^\el f$.  We consider the random field
$(X_\el = A^{\el} f, \, \el \in \Z^2)$, with $f \in AC_0(\T^\rho)$.

A sufficient condition for $f$ with 0 integral to be in $AC_0(\T^\rho)$ is $|\hat f(\k)| = O(\|k\|^{-\beta}), \text{ with } \beta > \rho$.

For a.e. $\omega$ the following asymptotic variance exists:
$$\lim_n \frac1{n\log n}\|\sum_{k=1}^{n}A^{Z_k(\omega)} f\|_2^2 = C_0 \, \sum_{\k \in \Z^d} \langle T^\k f \, f \rangle,$$
where the constant $C_0$ is defined in Subsection \ref{varRWSect1}.

The following quenched FCLT extends for the torus the CLT proved in \cite{CohCo17}. 
Remark that the CLT is proved therein for a general compact abelian group. 
The extension to a functional version of the CLT holds in this general case when $f$ is trigonometric polynomial.
\begin{thm} \label{rwFCLTB2} Let $(Z_n)$ be a 2-dimensional reduced centered random walk with a finite moment of order 2 
and let $f$ be a real function in $AC_0(\T^\rho)$ with spectral density $\varphi_f$ and $\varphi_f(\0) \not = 0$. 
Denoting by $S_{n}^\omega (f) :=\sum_{k=1}^{n}A^{Z_k(\omega)} f$ the sums along the r.w., the process
$\Big(\frac1{\sqrt{n\log n}}S_{\lfloor nt\rfloor}^\omega(f)\Big)_{t\in[0,1]}$ satisfies a FCLT for a.e. $\omega$.
\end{thm}
\proof {\it 1) Convergence of the finite dimensional distributions}

1a) First suppose that $f$ is a trigonometric polynomial: $f = \sum_{\k \in \Lambda} c_\k(f) \, \chi_\k$, where $(\chi_\k, \k \in \Lambda)$ 
is a finite set of characters on $\T^\rho$ and $\chi_{0}$ the trivial character. 

We use Proposition \ref{CltSum1}: (\ref{ortho1}) follows from (\ref{orth1}) and Lemma \ref{VarrwLem1}. For (\ref{negligCum3}), we have to show
\begin{eqnarray}
&C^{(r)}(\sum_{\el \in \Z^2} \, w_n(\omega, \el) \, T^\el f) = o((n \ln n)^{r/2}), \, \forall r \geq 3,  \text{ for a.e. } \omega. \label{cumulBdr0}
\end{eqnarray}
We apply Theorem \ref{Leonv0}. Let us check (\ref{smallCumul0}). For $r$ fixed, the function 
$(\n_1, ..., \n_r) \to m_f(\n_1, ..., \n_r) := \int_X \, T^{\n_1} f \cdots T^{\n_r} f \, d\mu$ takes a finite number of values, 
since $m_f$ is a sum with coefficients 0 or 1 of the products $c_{k_1} ... c_{k_r}$ with $k_j$ in a finite set. 
The cumulants of a given order take also a finite number of values according to (\ref{cumFormu1}).

Therefore, since mixing of all orders implies $\underset{\max_{i,j} \|\el_i - \el_j\| \to \infty} \lim \, C(T^{\el_1} f,..., T^{\el_r} f)  = 0$
by Proposition \ref{skTozeroLem}, there is $M_r$ such that $C(T^{\el_1} f,..., T^{\el_r} f) = 0$ if $\max_{i,j} \|\el_i - \el_j\| > M_r$.
The end of the proof is then like in Theorem \ref{FCLTindFil}.

1b) For $f \in AC_0(\T^\rho)$, using Proposition \ref{resum} and Lemma \ref{tclRegKer}, the convergence follows by approximation 
of $f$ by a squence of  trigonometric polynomials $f_L$ in such a way that $\lim_L \varphi_{f - f_L} (0)= 0$.

{\it 2) Moment of order 4 and tightness}

We use Proposition \ref{rwFCLTB}. 
Taking into account Remark \ref{remSum1}, it suffices for the tightness to take for $f$ a character and show that the bounds are independent of the character.

Let $\chi_v$ be a character on the torus $\T^\rho$, $\chi_v: x \to \exp(2 \pi i \langle v, x\rangle)$, where $v \in \Z^\rho \setminus \{0\}$.

For an interval $J =[b, b+k] \subset [1, n]$, we have: 
\begin{eqnarray*}
&\|\sum_{i= b }^{b+k} A^{Z_i} \chi_v \|_4^4 = 
\#\{(i_1, i_2, i_3, i_4) \in J^4: \ (A^{Z_{i_1}} -A^{Z_{i_2}}+A^{Z_{i_3}}-A^{Z_{i_4}}) \, v =0 \}.
\end{eqnarray*}
This number is less than $2G^2(\omega, b , k) + H(\omega, b , k)$, with
\begin{eqnarray*}
G(\omega, b , k) &&:= \#\{(i_1, i_2) \in J^2: \, (A^{Z_{i_1}} - A^{Z_{i_2}}) \v = 0 \},\\
H(\omega, b , k) &&:= \#\{(i_1, i_2, i_3, i_4) \in J^4: \, (A^{Z_{i_1}} - A^{Z_{i_2}} + A^{Z_{i_3}} - A^{Z_{i_4}}) \v= 0 \},
\end{eqnarray*}
where above in $H$ we count the number of solutions without vanishing proper sub-sums.

By assumption of total ergodicity,  if $(A^{Z_{i_1}} - A^{Z_{i_2}}) \v = 0$, then $Z_{i_1} = Z_{i_2}$, so that $G$ is the number of self-intersections 
of the r.w. on $[b, b+k]$: $G(\omega, b , k) = V(\omega, [b, b+k[) = V_k(\theta^b\omega)$.

For $H$, by Corollary \ref{nbsolu1}, there is a finite set $F$ (independent of $\v$) such that
$$H(\omega, b , k)  \leq \#\{(i_1, i_2, i_3, i_4) \in J^4: 
\ \bigl(Z_{i_1}(\omega) - Z_{i_4}(\omega), Z_{i_2}(\omega) - Z_{i_4}(\omega), Z_{i_3}(\omega) - Z_{i_4}(\omega)\bigr) \in F\}.$$
Therefore, with the notation (\ref{wn123}),  $H(\omega, b , k) \leq \sum_{(\el_1, \el_2, \el_3) \in F} W_k(\theta^b \omega, \el_1, \el_2, \el_3)$.

By (\ref{majnumb3}) in Lemma \ref{majwnm}, there exists a positive integrable function $C_3$ such that 
$W_n(\omega, \el_1, \el_2, \el_3) \leq C_3(\omega) \, n \, (\ln n)^5, \, \forall n \geq 1$, which implies
$H(\omega, b , k) \leq (\Card \, F) \ C_3(\theta^b \omega) \, k \, (\ln k)^5$.
Remark that the bounds do not depend on the character, but only on the matrices $A_1, A_2$.

Since $G$ and $H$ are super-additive (Condition (\ref{supadd1})), the tightness property follows now from Proposition \ref{rwFCLTB} with $\gamma = 5$. \eop

\begin{rems}
1) An analogous result is valid for any transient random walk in dimension $d \geq 1$, with the standrad normalisation by $\sqrt n$. 
In this case, if the observable is non null a.e., the asymptotic variance for the sums along the r.w. is different from 0 (cf. Subsection \ref{varRWSect1}). 

2) For automorphisms of a torus $\T^\rho$, in  the recurrent 2-dimensional model studied above, if $f$ satisfies the regularity condition 
$|c_f(\k)| = O(\|k\|^{-\beta})$, with $\beta > \rho$ , then the asymptotic variance is given by $\varphi_f(0)$ and it is null, 
if and only if $f$ is a mixed coboundary: there are continuous functions $u_1, u_2$ such that $f= (I- A_1) u_1 + (I- A_2) u_2$
(\cf \cite{CohCo17}).
\end{rems}

\vskip 3mm
\section{\bf Exponential mixing of all orders} \label{expMixSect}

\subsection{\bf FCLT and exponential mixing of all orders} \label{expMix}

\

Our last example is given by commuting translations on homogeneous spaces. It relies on recent results on the exponential mixing of all orders for flows
on homogeneous spaces shown in \cite{BjEiGo17} and their application to the CLT in \cite{BjGo20}. 
Closely following the latter reference, we recall first the notion of exponential mixing of all orders 

{\bf Exponential mixing of all orders}

Let $\cal H$ be a group with a left invariant distance $d$. Let $h \to T^h$ be a homomorphism of $\cal H$ in the group of a measure preserving invertible transformations  of 
a probability space $(E, \cal B, \mu)$.
We denote by ${\cal A}$ a sub-algebra in $L^\infty(E, \mu)$ which is $\cal H$-invariant.
Let ${\cal N} = (N_s)$ be a family of semi-norms on ${\cal A}$, indexed by positive integers $s$.

The following conditions are assumed to hold (with constant factors depending only on $s$, all denoted by $C_s$) for all $s > 1$ and all $f, g \in {\cal A}$:
\hfill \break 1) $N_s(f) \leq C_s N_{s+1}(f)$; 2) $\|f\|_{L^\infty} \leq C_s \, N_s(f)$; 3) $N_s(f g) \leq C_s \, N_{s+1}(f) \, N_{s+1}(g)$;
\hfill \break 4) there exists $\zeta_s > 0$ such that $N_s(T^h f) \leq C_s \, e^{\zeta_s d(h,e)} N_s(f), \, \forall h \in \cal H$.

Let $r > 2$ be an integer. For  $\underline H = (h_1, ... , h_r) \in \cal H^r$, we set $d_r( \underline H) := \min_{i \not =j} d(h_i, h_j)$. 

\begin{defi} \label{rexpMixDef} We say that the $\cal H$-action on $(E, \mu)$ is exponentially mixing of order $r$, with respect to $d$ and $({\cal A}, {\cal N})$, 
if there exist $\delta_r > 0$ and an integer $s_r > 0$ such that for all $s > s_r$ and $f_1, ... , f_r \in {\cal A}$,
\begin{eqnarray}
|\mu(\prod_{i=1}^r T^{h_i} f_i) - \prod_{i=1}^r \mu(f_i)| \leq C \,  e^{-\delta_r d_r(\underline H)} \prod_{i=1}^r N_s(f_i),  \label{rexpMix}
\end{eqnarray}
for all $\underline H = (h_1, ... , h_r) \in \cal H^r$. The constant $C$ depends only on $r$ and $s$.
\end{defi}
We may assume that $\zeta_s$ is increasing with $s$, $\delta_r$ decreasing with $r$ and $\delta_r < r \zeta_s, \, \forall r, s$.

In what follows, we will consider, for $d \geq 2$, a measure preserving $\Z^d$-action on a probability space $(E, \cal B, \mu)$ 
generated by $d$ commuting invertible maps $T_1, ..., T_d$. Therefore the group $\cal H$ in Definition \ref{rexpMixDef} is going to be the group $\Z^d$ 
still denoted also by $\cal H$.

If $(Z_n)$ is a random walk on $\Z^d$, then we get a random walk $(T^{Z_n(\omega)})$ on the group of measure preserving invertible transformations on $(E, \cal B, \mu)$.

With a distance on $\Z^d$ associated to a norm equivalent to the Euclidean norm, the volume of a big ball is of order the number of integral points in the ball.
It is important to relate this distance to the distance $d$ of Definition \ref{rexpMixDef}. We will assume that the action satisfies:
\begin{hyp} \label{ZdMix} For an $\cal H$-invariant sub-algebra ${\cal A}$ in $L^\infty(E, \mu)$ and a family ${\cal N} = (N_s)$ of semi-norms on ${\cal A}$,
we assume that the $\Z^d$-action on $(E, \cal B, \mu)$ is exponentially mixing of order $r$ for every $r \geq 2$ in the sense of Definition \ref{rexpMixDef}
with a distance $d$ equivalent to the Euclidean distance.
\end{hyp}
With this assumption, for simplicity of notation we can assume that the distance $d$ in Definition \ref{rexpMixDef} applied to $\Z^d$ is the Euclidean distance on $\Z^d$.

{\bf Spectral density and cumulants}

Let $f$ be a centered function in $\cal A$. Its spectral density is $\varphi_f = \sum_{\el \in \Z^2} a_\el e^{2\pi i \langle\el, \, \t\rangle}$,
with $a_\el = \langle T^\el f, f\rangle$. The absolute summability $\sum_\el |a_\el| < \infty$ is a consequence of (\ref{rexpMix}) for $r=2$.

Let $J = [b, b+k] \subset [1, n]$. With the notation of \ref{selfInter}, we have
\begin{eqnarray*}
\int |\sum_{j \in J} T^{Z_j(\omega)} f|^2 \, d\mu &=& \int_{\T^d} |\sum_{j \in J} e^{2\pi i \langle Z_j(\omega), \,\t\rangle}|^2 \, \varphi_f(\t) \, d\t
 \leq \sum_{\el \in \Z^2} |a_\el| \, V(\omega, J, \el) =: G(\omega, b, k). 
\end{eqnarray*}
For a fixed $\omega$, the bound $ G(\omega, b, k)$ is super-additive. By (\ref{moment4}), we have
\begin{eqnarray*}
&&\E_\mu(|\sum_{j \in J} T^{Z_j(\omega)} f|^4) = 3 \bigl(\int (\sum_{j \in J} T^{Z_j(\omega)} f)^2\, d\mu \bigr)^2  + C^{(4)}(\sum_{j \in J} T^{Z_j(\omega)} f).
\end{eqnarray*}
The next proposition shows that the cumulants $C(T^{h_1} f,..., T^{h_r} f)$ are small for all well-separated $r$-tuples $\underline H = (h_1, . . . , h_r)$.
The constants $\delta_r$ and $\zeta_r$ are those of the beginning of the subsection.
\begin{prop} \label{majDeltaQ} (cf. Proposition 6.1 in \cite{BjGo20}) For $r \geq 2$, let $Q$ be a partition of $\{1, 2, ..., r\}$ with $|Q| \geq 2$ and
let $s$ be an integer $> s_r + r$. Then, if $\underline H \in \Delta_Q(\alpha, \beta)$ with $0 < \alpha < \beta$, for $f$ centered in $\Cal A$ we have
\begin{eqnarray}
|C(T^{h_1} f,..., T^{h_r} f)| \leq C e^{- (\delta_r \, \beta - r \, \alpha \, \zeta_s)} \, N_{s}(f)^r.
\end{eqnarray}
where the constant $C$ depends only on $r$ and $s$.
\end{prop}
\proof Let us give the proof for $r=4$ and in the case of the configurations $\underline H$ of type 
$S_3 = \bullet \bullet \ \ \  \bullet \ \ \ \bullet$, that is (up to a permutation): $\{\{h_i, h_k\}, \{ h_\ell\}, \{  h_r\}\}$ 
with, for some $j$, $d(h_i, h_k) < 3\beta_j$ and $d(h_\ell, \{h_i, h_k\}) > \beta_{j+1}$, $d(h_r, \{h_i, h_k\}) > \beta_{j+1}$,  $d(h_\ell, h_r) > \beta_{j+1}$.
 
We may write the formula for the cumulants in the following way:
\begin{eqnarray*}
&& C(T^{h_i} f, \,T^{h_k} f, \,T^{h_\ell}f, \, T^{h_r} f)  = (A) -  (B) - (C) - (D), \text{ with }\\
&&(A) = \E \bigl(T^{h_i}  f \, T^{h_k} f \, T^{h_\ell}  f \, T^{h_r}  f\bigr), \, (B) =  \E (T^{h_i}  f \, T^{h_k} f) \, \E(T^{h_\ell}  f \, T^{h_r}  f),\\
&&(C) = \E(T^{h_i}  f \, T^{h_\ell} f) \, \E(T^{h_k}  f \, T^{h_r}  f), \, (D) = \E(T^{h_i}  f \, T^{h_r} f) \, \E(T^{h_k}  f \, T^{h_\ell}  f).
\end{eqnarray*}
We use the exponential mixing of order 2 for $(B), (C), (D)$ and of order 3 for $(A)$. More precisely, we have 
$(A) = \E \bigl((f \, T^{h_k - h_i} f) \, T^{h_\ell  - h_i} f \, T^{h_r  - h_i} f\bigr)$ and
$$N_s(f \, T^{h_k - h_i} f) \leq C N_{s+1}(f) \, N_{s+1}(T^{h_k - h_i} f) \leq C e^{3 \zeta_s \, \beta_j} N_{s+1}(f) \, N_{s+1}(f).$$
Therefore, by the 3-mixing applied to $f \, T^{h_k - h_i} f, \, T^{h_\ell  - h_i} f, \, T^{h_r  - h_i} f$, we have
$$(A) \leq  C \, e^{3\zeta_s \, \beta_j} \, e^{-\delta_r \,\beta_{j+1}} \, N_{s+1}(f)^2 \,  N_{s}(f)^2. \eop$$

Let us now fix $\gamma > 0$. We define $\beta_{j}$ by $\beta_0 = 0$ and 
$\beta_{j+1} = 3 r \beta_{j} \zeta_s / \delta_r + \gamma$, for $j > 0$.

As $\delta_r < r \zeta_s$, we have $\beta_{j+1} > 3 \beta_{j}$. Moreover $\beta_{r} = \gamma \sum_{j=0}^{r-1} (3 r \, \zeta_s/ \delta_r)^j =: \gamma c_{r,s}$.

Let $\nu = \sum_\el \nu(\el) \delta_\el$ be a positive finite measure on $\cal H= \Z^d$.
The following bound, where the first term comes from the clustered configurations and the second term from the well separated configurations,
results from Proposition \ref{majDeltaQ} and from (\ref{unionDelta}):
\begin{prop} (\cf Proposition 5.2 in \cite{BjGo20}) For every $r \geq 2$, there exist  $C$ and $c= c_{r,s} > 0$ such that,
for all $\gamma > 0$ and $f \in \Cal A$,
\begin{eqnarray}
|C^{(r)}(\nu*f)| \leq C \bigl(\int_{\cal H} \nu(B(h, c \gamma))^{r-1} \, d\nu(h) + e^{-\delta_r \gamma} \, \|\nu\|^r \bigr) \, N_s(f)^r. \label{majcumr}
\end{eqnarray}
\end{prop}
Under the assumption of exponential mixing of all orders, it is shown in \cite{BjGo20} that the CLT holds for $\nu_n*f$, when $f \in \Cal A$
and $(\nu_n)$ a sequence of measures satisfying a certain condition.

In our framework, $\nu$ is the measure $\nu_n^\omega = \sum_{j=0 }^{n-1} \delta_{Z_j(\omega)}$, where $(Z_j)_{j \geq 0}$ is a r.w. on $\Z^2$. 
Its mass is $\|\nu_n^\omega\| = n$.
The convolution ``$\nu*f$'' in (\ref{majcumr}) for $\nu_n^\omega$ means ``$\sum_{j = 0}^{n-1} f(T^{Z_j(\omega)} .)$''.
 
Let us assume that the r.w. $Z$ is a centered random walk on $\Z^2$ with moments of all orders.
By (\ref{DPRZbound}) (\cf \ref{selfInter}), it holds with a constant $K(\omega)$ finite for a.e. $\omega$, uniformly in $h \in \Cal H$:
\begin{eqnarray*}
&\nu_n^\omega(B(h,c\gamma))=\#\{j\le n: Z_j(\omega)\in B(h, c\gamma)\}  \leq \sup_\el w_n(\omega, \el) \, \Card(B(h, c\gamma)) \leq K(\omega) (c\gamma)^2\log^2 n.
\end{eqnarray*} 
Due to Hypothesis \ref{ZdMix} and by (\ref{majcumr}), this gives with a constant $C_1$ depending only on $r$ and $s$, 
\begin{eqnarray*}
|C^{(r)}(\nu_n^\omega*f)| \leq C_1 \bigl( K(\omega)^{r-1} \, (\ln n)^{2(r-1)} (c\gamma)^{2(r-1)} n+ e^{-\delta_r \gamma} \, n^r \bigr) \, N_s(f)^r.
\end{eqnarray*}
Taking $\gamma ={r-1 \over \delta_r} \, \ln n$, the first term above gives the bound $(c {r-1 \over \delta_r})^{2(r-1)} \, (\ln n)^{4(r-1)} n$ and the second term is $n$.
It follows, for a constant $C_2$ depending on $r$ and $s$:
\begin{eqnarray*}
|C^{(r)}(\nu_n^\omega*f)| \leq C_2 \, N_s(f)^r \, K(\omega)^{r-1}  \, (\ln n)^{4(r-1)} \, n.
\end{eqnarray*}
Likewise for an interval $J = [b, b+k]$ and $\nu_J^\omega = \sum_{j \in J } \delta_{Z_j(\omega)}$, we get
\begin{eqnarray*}
|C^{(r)}(\nu_J^\omega*f)| \leq  C_2 \, N_s(f)^r \, K(\theta^b \omega)^{r-1}  \, (\ln k)^{4(r-1)} \, k.
\end{eqnarray*}
So we get the same type of upper bound for the cumulants as for the automorphisms of the torus. 
Therefore we have convergence of the finite dimensional distributions and tightness.

For a $\Z^2$-action satisfying the exponential mixing condition on an algebra of functions as presented at the beginning of this section, 
we can state now a functional version of a CLT result for the summation along a random walk.
The result is formulated for $\Z^2$, but an analogous result can be proved with the same method for a transient random walk in dimension $\geq 1$.

Let $(Z_n)$ be a 2-dimensional aperiodic centered random walk with moments of all orders. Then we have:
\begin{thm} \label{rwFCLTExp} Let $f$ be a real centered function in $\cal A$ with spectral density $\varphi_f$ such that $\varphi_f(\0) \not = 0$. 
Under Hypothesis \ref{ZdMix}, denoting by $S_{n}^\omega (f) :=\sum_{k=1}^{n}T^{Z_k(\omega)} f$ the sums along the r.w., 
the process$\Big(\frac1{\sqrt{n\log n}}S_{\lfloor nt\rfloor}^\omega(f)\Big)_{t\in[0,1]}$ satisfies a FCLT for a.e. $\omega$.
\end{thm}

\vskip 3mm
\subsection{\bf Translations on homogeneous spaces} \label{homog}

\ 

The following example is an action which is exponentially mixing of all orders on an algebra $\cal A$ of functions according to \cite{BjEiGo17}.

We take the group $G = \rm{SL}(n, \R)$ and a lattice $\Gamma$ in $G$, i.e.,  a discrete subgroup such that  $G /\Gamma$ has a finite volume 
for the measure $\mu$ induced by the Haar measure of $G$, for example $\rm{SL}(n, \Z)$.  The space $E$ is the quotient $G /\Gamma$. 
The action on $E$ will be given by left multiplication $g \Gamma \to hg\Gamma$ where $h$ is in the diagonal subgroup $D$ of $G$.

The algebra $\cal A$ in the example is the algebra of $C^\infty$-functions with compact support on $G /\Gamma$, and $\Cal N$ is a family of Sobolev norms as in \cite{BjEiGo17}.

Observe that it is not true that in general for an invariant distance on an abelian or a nilpotent group, the growth of the balls is sub-exponential. 
(Notice that Theorem 6.8.1 in \cite{CeCo10} as referred in \cite{BjGo20} is shown for the growth corresponding to the length associated to a set of generators, 
and not for the growth of the balls associated to any invariant distance). 
Therefore in the example, we need to explicit the distance .

{\it Left invariant pseudo-metric on $G = \rm{SL}(n, \R)$}

Recall that there is a distance on $D$, induced by a canonical left invariant distance on $G$.
The left invariant distance $d_G$ defined on $G$ is comparable on $D$ to the `pseudo-metric' $\delta(., .)$ defined as follows:

For $A \in G$, let $|||A|||= \max_{x \not = 0} { \|A x\| \over \|x\|}$ be the norm of the $n \times n$ matrix $A$ as operator on $\R^n$ endowed with the euclidian norm. 
Since the determinant is 1, $A$ has an eigenvalue of modulus $\geq 1$, which implies $|||A||| \geq 1$. For $A, B \in G$, we put $\delta (A, B) = \log (|||A^{-1} \, B|||)$. 

Clearly, $\delta (A, B)  \geq 0$, the triangular inequality is satisfied by sub-multiplicativity of the operator norm, and $\delta$ is left invariant on $G$. 
If $|||A|||= 1$, the iterates of $A$ are bounded, so all eigenvalues must have a modulus $\leq 1$. As the determinant is 1, the modulus of the eigenvalues is 1.
Now considering the Jordan form of $A$ over $\C$, it must be diagonal. Finally we conclude that $\{A: |||A||| =1\}$ is the orthogonal group in $G$. 

We take two elements $A_1 = \exp(U_1)$, $A_2= \exp(U_2)$, where $U_1, U_2$, in the sub-algebra $\mathfrak D$
corresponding to $D$ in the Lie algebra $\mathfrak G$ of $G$, are such that $U_1, U_2$ generate a 2-dimensional vectorial space of $\mathfrak D$. 
The group $(A_1^{\ell_1} A_2^{\ell_2}, \, (\ell_1,\ell_2) \in \Z^2)$ yields a totally ergodic action on $G / \Gamma$. 

For instance in $\rm{SL}(3, \R)$, we can take with $a_1, a_2 > 1$:
$A_1 = \left(
\begin{matrix} a_1 & 0 & 0 \cr 0 & 1 & 0 \cr 0 & 0 & a_1^{-1} \cr
\end{matrix} \right), \ \ A_2  = \left(
\begin{matrix} 1 & 0 & 0 \cr 0 & a_2 & 0 \cr 0 & 0 & a_2^{-1}\cr
\end{matrix}
\right)$.

The distance $d_G(A_1^{\ell_1} A_2^{\ell_2}, Id)$ is equivalent to $\|\el\| = (\ell_1^2 + \ell_2^2)^\frac12$.
The measure $\nu_n(B(Id, R))$ is the counting measure with some weight applied to the ball, 
therefore up to the weight it is the number of elements of the form $h_1^{\ell_1} h_2^{\ell_2}$ in the ball,
and finally the (weighted) number of integers $\el = (\ell_1, \ell_2)$ of norm $\leq R$.

By what precedes, Hypothesis \ref{ZdMix} is satisfied and Theorem \ref{rwFCLTExp} yields a functional CLT in the class of centered 
compactly supported $C^\infty$-functions for the action of a 2-dimensional random walk on the diagonal subgroup on $G/ \Gamma$.

An result analogous to Theorem \ref{rwFCLTExp} holds for the sums along a transient random walk: the only change is the estimate of the number of self-intersection,
(normalisation by $\sqrt n$).

In Theorem \ref{rwFCLTExp}, in particular in the example provided by homogeneous spaces, as in the CLT in \cite{BjGo20}, 
the statement says nothing about the non-nullity of the variance, for a 2-dimensional recurrent r.w.
The question of degeneracy of the asymptotic variance (for a recurrent r.w.) is the same as for the sums over squares: it depends of the nullity of $\varphi(\0)$.
This contrasts with the action of commuting automorphisms of a torus, for which we have seen that there is a description of the degenerate case in terms of mixed coboundaries.

In the case of a transient random walk on $\Z^d$, as noticed in Subsection \ref{varRWSect1}, if the observable is non null a.e., 
the asymptotic variance for the sums along a transient r.w. is different from 0.

\end{document}